\documentclass[pre,letterpaper,twocolumn,nofootinbib,showkeys]{revtex4}

\usepackage{amsmath}
\usepackage{graphicx}
\usepackage{hyperref}
\usepackage{listings}
\usepackage{algorithmic}
\usepackage{algorithm}
\usepackage{textcomp}
\usepackage{verbatim}
\usepackage{color}
\usepackage{multirow}

\date{\today}

\begin{document}

\title{An Optimized Sparse Approximate Matrix Multiply for Matrices with
Decay}

\preprint{LA-UR 11-06091}

\author{Nicolas Bock}
\email{nbock@lanl.gov}
\affiliation{Group T-1, Theoretical Division, Los Alamos National Laboratory,
Los Alamos, NM 87544}

\author{Matt Challacombe}
\email{mchalla@lanl.gov}
\affiliation{Group T-1, Theoretical Division, Los Alamos National Laboratory,
Los Alamos, NM 87544}

\newcommand{\SpAMM}{{\tt SpAMM}}

\keywords{Sparse Approximate Matrix Multiply; Sparse Linear Algebra; SpAMM;
Reduced Complexity Algorithm; Quantum Chemistry; $N$-Body; Matrices with Decay}


\begin{abstract}
We present an optimized single-precision implementation of the Sparse
Approximate Matrix Multiply (\SpAMM{}) [M.~Challacombe and N.~Bock, arXiv {\bf
1011.3534} (2010)], a fast algorithm for matrix-matrix multiplication for
matrices with decay that achieves an $\mathcal{O} \left( n \log n \right)$
computational complexity with respect to matrix dimension $n$.  We find that
the max norm of the error achieved with a \SpAMM{} tolerance below $2 \times
10^{-8}$ is lower than that of the single-precision {\tt SGEMM} for dense
quantum chemical matrices, while outperforming {\tt SGEMM} with a cross-over
already for small matrices ($n \sim 1000$).  Relative to naive implementations
of \SpAMM{} using Intel's Math Kernel Library ({\tt MKL}) or AMD's Core Math
Library ({\tt ACML}), our optimized version is found to be significantly
faster. Detailed performance comparisons are made for quantum chemical
matrices with differently structured sub-blocks. Finally, we discuss the
potential of improved hardware prefetch to yield 2--3x speedups.
\end{abstract}

\maketitle


\section{Introduction}

For large dense linear algebra problems, the computational advantage offered by
fast matrix-matrix multiplication can be substantial, even with seemingly small
gains in asymptotic complexity\footnote{Throughout, we refer to the asymptotic
bounds of an algorithm in big-$\mathcal{O}$ notation when using the term
``complexity'', while the term ``performance'' is reserved for the actual
runtime or cycle count of an implementation. Since performance is affected by
the complexity of the underlying algorithm \emph{and} computer and
implementation specific details, complexity and performance are not
necessarily the same.}.  Relative to conventional multiplication which is
$\mathcal{O} \left( n^{3} \right)$, Strassen's algorithm
\cite{springerlink:10.1007/BF02165411} achieves $\mathcal{O} \left( n^{2.8}
\right)$, while the Coppersmith and Winograd method \cite{Coppersmith1990251}
and the method of Williams' \cite{Williams:2012:MMF:2213977.2214056} are bound
by $\mathcal{O} \left( n^{2.3755} \right)$ and $\mathcal{O} \left( n^{2.3727}
\right)$, respectively. For these dense methods, balancing the trade-off
between cost, complexity and error is an active area of research
\cite{Demmel:1992:FastMM, Demmel:2007:FastMM, Yuster:2005:FastMM}. On the
other hand, large sparse problems are typically handled with conventional
sparse matrix techniques, with only small concessions between multiplication
algorithms.

\begin{figure}
\includegraphics[width=1.0\columnwidth]{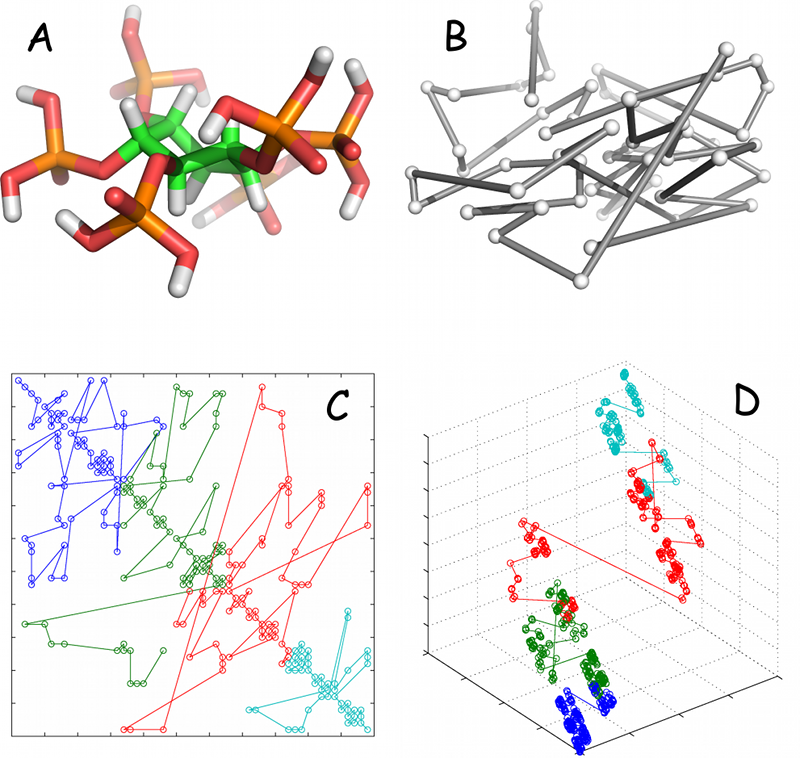}
\caption{\label{fig:SFC} Space filling curves (SFC) map atoms in Cartesian
space (A) onto an index that is locality preserving (B), leading to clustering
of matrix elements due to decay and the underlying physical system, shown for
the sparsified case (C). In addition the locality preserving mechanisms of the
SFC can be extended to 2-D ($ij$) data decomposition of matrices, (C), as well
as three-dimensional decomposition ($ijk$) in product space of the
matrix-matrix multiply (D).}
\end{figure}

\begin{figure}
\includegraphics[width=1.0\columnwidth]{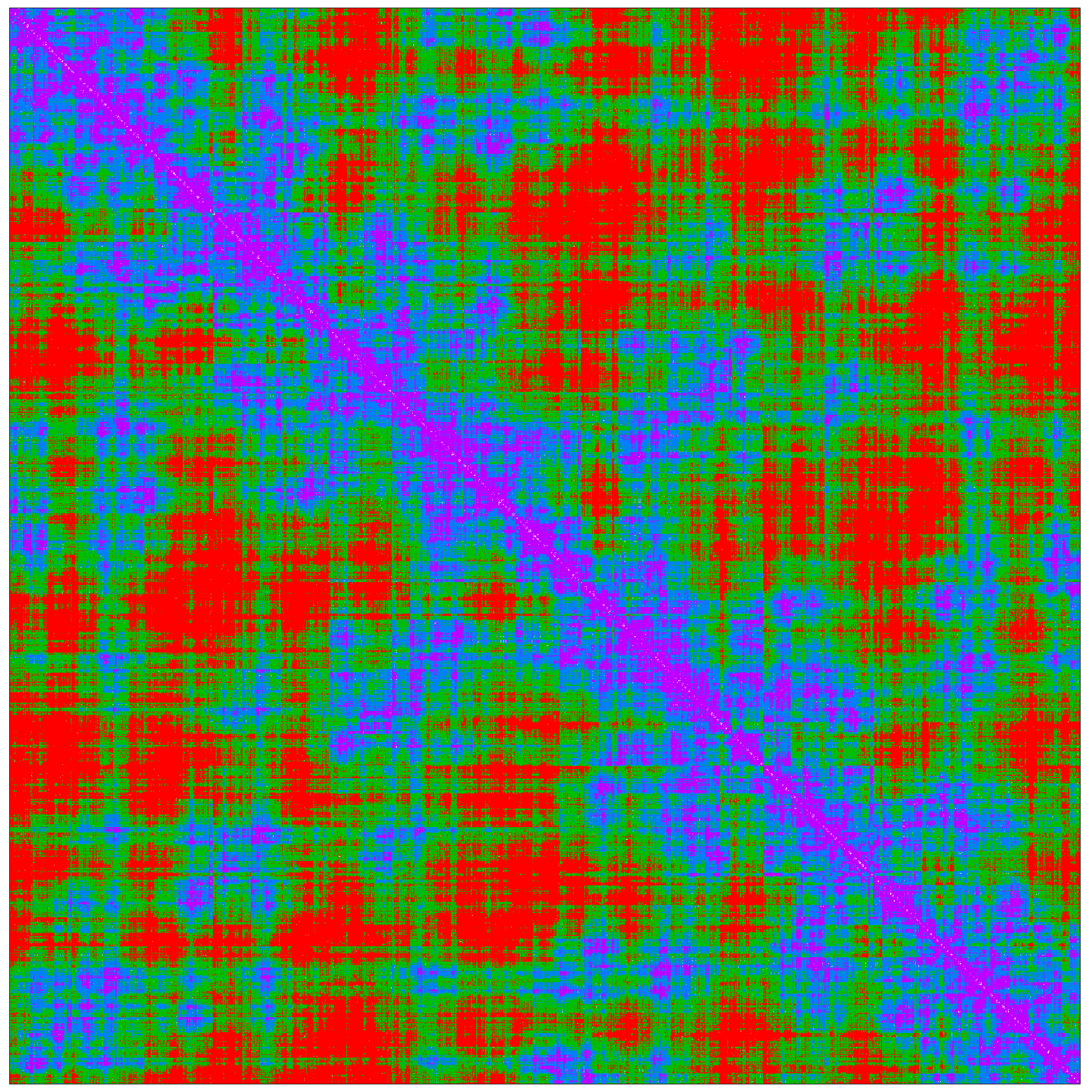}
\caption{\label{fig:sparsity_pattern} The decay of matrix element magnitudes of
a converged spectral projector (density matrix) for a $\left( \mathrm{H}_{2}
\mathrm{O} \right)_{300}$ water cluster at the RHF/6-31G${}^{**}$ level of
theory ($n = 7500)$, where the molecular geometry has been reordered with a
space filling Hilbert curve. The different colors indicate different matrix
element magnitudes; red: $\left[ 0, 10^{-8} \right)$; green: $\left[ 10^{-8},
10^{-6} \right)$; blue: $\left[ 10^{-6}, 10^{-2} \right)$; violet: $\left[
10^{-2}, 1 \right]$, corresponding to approximate exponential decay.}
\end{figure}

Intermediate to these regimes, a wide class of problems exist that involve
matrices with decay, occurring in the construction of matrix functions
\cite{benzigolub}, notably the matrix inverse \cite{demkomosssmith,
Benzi:2000:Inv}, the matrix exponential \cite{Iserles:2000:Decay}, and in the
case of electronic structure theory, the Heaviside step function (spectral
projector) \cite{McWeeny12061956, PhysRevB.58.12704, Challacombe:1999:DMM,
Challacombe:2000:SpMM, Benzi:2007:Decay, BenziDecay2012}. A matrix $A$ is said
to decay when its matrix elements decrease exponentially, as $\left| a_{i, j}
\right| < c \,\, \lambda^{\left| i-j \right|}$, or algebraically as $\left|
a_{i, j} \right| < c / ( \left| i-j \right|^{\lambda}+1 )$ with separation
$|i-j|$. In non-synthetic cases, the separation $|i-j|$ may correspond to an
underlying physical distance $|\vec{r}_{i}-\vec{r}_{j}|$, \emph{e.g.} of basis
functions, finite elements, \emph{etc.} \cite{BenziDecay2012}. Note, matrices
with decay are typically not sparse, but exhibit a structure that can be
exploited to yield an approximate sparse algorithm that exhibits a reduced
complexity.

Reordering the physical indices $\vec{r}_{i}$ along a locality preserving
space-filling curve (SFC), such as the Hilbert \cite{HilberCurve} or the Peano
curve \cite{peano1890courbe}, maps elements that are close in three dimensions
onto an ordered list \cite{Warren:1992:HOT, Warren:1995:HOTb,
Samet:2006:DBDS}, leading to effective clustering of submatrices by order of
magnitude \cite{Challacombe:2000:SpMM}.  Figure \ref{fig:SFC} illustrates this
point for quantum chemistry applications where the atom indices of the
molecule shown in panel A are ordered along the SFC, shown in panel B, leading
to the clustered structure of matrix element magnitudes shown schematically in
panel C of Fig.~\ref{fig:SFC} and specifically in
Fig.~\ref{fig:sparsity_pattern}.  Extending SFC ordering to the two
dimensional vector space of matrix elements $A_{ij}$, shown in panel C, and
the three dimensional product space of product contributions $A_{ik} B_{kj}$,
shown in panel D, extends the locality preserving properties of SFC to matrix
element storage and matrix product execution order. Historically, complexity
reduction for matrices with decay has been achieved through ``sparsification''
or dropping of small matrix elements and the use of conventional sparse matrix
kernels \cite{Galli:1996:ONREV, Goedecker98, Goedecker:1999:ONREV,
0953-8984-14-11-303, Goedecker:2003:ONREV, Li:2005:ONREV, bowler2011:ONREV}.
Alternatives such as divide and conquer exist \cite{yang1991direct,
yang1995density, lee1998linear, Nakano2007642}, where block diagonal
subsystems are ``glued'' together, but these methods are unable to retain long
range coherence effects that can manifest in the spectral projector, show for
example in Fig.~\ref{fig:sparsity_pattern} as \emph{cross-diagonal} bands.
Recently, the authors introduced a fast method for calculating the product of
matrices with decay, the Sparse Approximate Matrix Multiply (\SpAMM{})
\cite{ChallacombeBock2010}, which is different from truncation in the vector
space, ($ij$) and ($jk$), in that it hierarchically eliminates insignificant
contributions in the product space ($ijk$), leading to complexity reduction
through a sparse product space.  The recursive approach is similar to
Strassen's algorithm \cite{springerlink:10.1007/BF02165411} and also to
methods of improving locality across memory hierarchies pioneered by Wise and
coworkers \cite{Frens:1997:AMT:263767.263789, Wise:2001:LSM:568014.379559,
Lorton:2006:ABL:1166133.1166134, Wise1990282,
springerlink:10.1007/3-540-51084-2_9, Wise:1984:RMQ:1089389.1089398,
Wise:Ahnentafel}.

\SpAMM{} belongs to a broad class of generalized $N$-Body solvers, ``fast''
algorithms that achieve reduced complexity through hierarchical approximation,
achieving increasingly broad applicability and high performance
implementations with surprising commonalities. Related solvers include the
database join operation \cite{Mishra:1992:JPR:128762.128764,
Schneider:1990:TPC:94362.94514, Chen:2007:IHJ:1272743.1272747,
Kim:2009:SVH:1687553.1687564}, gravitational force summation
\cite{10.1109/SUPERC.1992.236647, Warren:1993:PHO:169627.169640,
Warren:1995:HOTb}, the Coulomb interaction \cite{Challacombe:1996:QCTCb,
schwegler:8764}, the fast Gau\ss{} transform \cite{greengard:79, strain:1131,
baxter:257, yang2003:FGT, Wan20067}, adaptive mesh refinement
\cite{Berger1984484, bell:127, Berger198964}, the fast multipole method
\cite{Greengard1987325, Greengard1997, Cheng1999468,
springerlink:10.1007/BFb0089775} and visualization applications
\cite{CGF:CGF1554, CGF:CGF1775, Lefebvre:2006:PSH:1141911.1141926,
Lefebvre:2006:PSH:1179352.1141926, AVRIL:2009:HAL-00412870:1, Zou2008,
694268}.

\SpAMM{} is similar to the $\mathcal{H}$-matrix algebra of Hackbusch and
co-workers \cite{Hackbusch2002, Grasedyck2003} and the Hierarchically
Semiseparable (HSS) representation of Chandrasekaran, Gu and co-workers
\cite{Chandrasekaran2005, Chandrasekaran2006}, where off-diagonal sub-matrices
are treated hierarchically as reduced rank factorizations (truncated SVD),
typically structured and grouped to reflect properties of the underlying
operators. For problems with rapid decay, truncated SVD may behave in a
similar way to simple dropping schemes. However, \SpAMM{} is different than
the $\mathcal{H}$-matrix or the HSS algebra as it achieves separation uniquely
in the product space and does not rely on a reduced complexity representation
of matrices.  For very slow decay, where \SpAMM{} is ineffective, the
$\mathcal{H}$-matrix or HSS algebra may certainly offer an alternative.

In Ref.~\cite{ChallacombeBock2010}, a na\"{\i}ve implementation of \SpAMM{}
was developed using recursion and a simple unoptimized $4 \times 4 \times 4$
multiplication at the lowest level. In this article we describe an optimized,
single-precision, non-recursive kernel of \SpAMM{} that can be introduced at
any level of a recursive scheme targeting full precision of the conventional
{\tt SGEMM} kernel. Precision of the {\tt SGEMM} corresponds roughly to many
approximate double-precision $\mathcal{O} (n)$ schemes, which yield $\approx
7$ digits in the total energy and $\approx 5$ digits in the forces, with an
additional two-fold savings in time and space.  The optimized algorithm
includes a hand-coded assembly kernel using single-precision SSE instructions
for $4 \times 4 \times 4$ matrix products, demonstrating significantly higher
performance than any na\"{\i}ve implementation using vendor tuned {\tt BLAS}
kernels, while yielding superior error control.

Our article is organized as follows: In Sec.~\ref{sec:algorithm_and_theory} we
describe briefly the \SpAMM{} algorithm for the matrix multiply. In
Sec.~\ref{sec:implementation}, we describe in detail our high-performance
implementation. In Sec.~\ref{sec:methodology}, we describe the benchmark setup
and methodologies used. In Sec.~\ref{sec:results}, we present our results for
benchmarks on representative density matrices from quantum chemistry
calculations performed on two hardware platforms, Intel Xeon X5650 and AMD
Opteron 6168.  In Sec.~\ref{sec:discussion}, we discuss the results and,
finally, in Sec.~\ref{sec:conclusions}, we present our conclusions.

\section{Algorithm and Theory}
\label{sec:algorithm_and_theory}

The {\tt SGEMM} (single precision) function of {\tt BLAS} calculates the
following expression,
\begin{equation}
C = \alpha A B + \beta C,
\end{equation}
where $A$, $B$, and $C$ are $m \times k$, $k \times n$, and $m \times n$
matrices, respectively, and $\alpha$ and $\beta$ are scalar parameters.
\SpAMM{} calculates the same but is optimized for matrices with decay. In
\SpAMM{} we store an $m \times n$ matrix in a quaternary tree (quadtree) data
structure \cite{Finkel1974}, an idea developed by Wise \emph{et al.}
\cite{springerlink:10.1007/3-540-51084-2_9, Wise:1984:RMQ:1089389.1089398,
Wise:Ahnentafel}, Beckman \cite{beckman1993parallel}, and Samet
\cite{Samet:1990:DAS:77589, Samet:2006:DBDS}.  The quadtree resolves the
matrix decay pattern hierarchically and is particularly well suited for the
representation of matrices with a clustered structure such as that shown in
Fig.~\ref{fig:sparsity_pattern} and explained in Fig.~\ref{fig:SFC}.
Furthermore, computational cost for quadtree access exhibits favorable scaling
with respect to matrix size \cite{Wise1990282}.

Compared to the more traditional Compressed Sparse Row or Column (CSR or CSC,
respectively) data structure, which is typically stored in contiguous arrays,
the quadtree can easily be allocated non-contiguously and allows for insertion
and deletion without requiring data copy operations.  Note that even in the
case of more advanced CSR data structures the pointer table is allocated
contiguously, see for instance the BCSR data structure of Challacombe
\cite{Challacombe:2000:SpMM} which uses blocked data, or the doubly compressed
sparse row format (DCSR) of Bulu\c{c} \emph{et al.} \cite{Buluc2008} developed
for the efficient representation of hyper-sparse matrices.

In addition, CSR based implementations face challenges in parallel, currently
significantly suppressing performance scaling to large processor counts. In
the context of electronic structure calculations for instance, it has been
observed that data movement in parallel CSR summation can dominate the
parallel construction of sparse Hamiltonians \cite{10.1063/1.1568734}.
Similarly, the parallel sparse matrix-matrix multiply ({\tt SpGEMM}) of
Bulu\c{c} \emph{et al.} \cite{Buluc:2008:SpMM, buluc2011parallel} becomes
communication bound with increasing processor count due to the 2-D
decomposition and {\em random permutation} of matrix rows and columns that
simplify communication and load balancing, but lead to communication cost
scaling as in the dense (Cannon's algorithm \cite{cannon1969cellular} or SUMMA
\cite{van1997summa}) case, $\mathcal{O} \left( \alpha \sqrt{p} + \beta \, c \,
n / \sqrt{p} \right)$, where $\alpha$ represents communication latency,
$\beta$ communication bandwidth, and $c$ the average number of non-zero
elements per row.

There is also a fundamental philosophical difference between \SpAMM{} and CSR
based matrix-matrix multiply algorithms with {\tt SpGEMM} as a state of the
art example. On the one hand, \SpAMM{} exploits \emph{locality} manifest in
clustering of matrix elements by order of magnitude
(Fig.~\ref{fig:sparsity_pattern}) to achieve a reduced complexity through
hierarchical approximation in the product space. On the other hand, the
parallel {\tt SpGEMM} exploits \emph{randomization} and corresponding loss of
locality to achieve load balance. While Bulu\c{c}'s idea has been deployed by
a large scale quantum chemistry program \cite{VandeVondele2012}, the impact of
data reordering and loss of locality on scalability would appear to also
suffer from communication bottlenecks with same formal complexities as
Bulu\c{c}'s.  It remains for us to demonstrate that these complexities can be
overcome in the \SpAMM{} approach by retaining locality and achieving load
balance instead in overdecomposition of recursive task space.  In addition,
the {\tt SpGEMM} approach would strongly suppress performance of quantum
chemical algorithms due to total loss of (matrix) locality,  as well as
interoperability problems due to global ordering and reordering of data.

The non-contiguous allocation of a quadtree and the SFC ordering of matrix
elements should therefore be more suited for integration with other quantum
chemical solvers and in parallel implementations where data and work
distribution in a non-shared memory environment is important and data
repacking undesirable.

\subsection{Quadtree Data Structure}
\label{sec:quadtree}

Given an $m \times n$ matrix, the depth of the quadtree is given by
\begin{equation}
d = \lceil \max \left( \frac{\log \left( m / n_{b} \right) }{\log 2},
\frac{\log \left( n / n_{b} \right) }{\log 2} \right) \rceil,
\end{equation}
where $n_{b} \times n_{b}$ is the size of the dense submatrices stored at the
bottom of the tree, chosen for performance.  To simplify recursion, the matrix
is zero padded to square shape of size $n_{p} = n_{b} \,\, 2^{d}$. The
computational cost of the lookup, insertion, and deletion of matrix elements
is $\mathcal{O} \left( \log n_{p} \right)$, independent of the decay pattern.
The root tier is denoted as $t = 0$ and at each tier $t < d$ of the quadtree,
a node links to the subnodes of the $2 \times 2$ submatrix underneath it,
\begin{equation}
A^{t} = \left(
\begin{array}{cc}
A^{t+1}_{11} & A^{t+1}_{12} \\
A^{t+1}_{21} & A^{t+1}_{22}
\end{array}
\right).
\end{equation}
If a submatrix is zero, its link is set to {\tt NULL} and the tree terminates
at that point.  Finally, at the bottom tier $t = d$, each node stores a dense
$n_{b} \times n_{b}$ submatrix.

\subsection{\SpAMM{}}
\label{sec:spamm}

In each node we also store the Frobenius norm of its submatrix, which is given
by
\begin{equation}
\lVert A^{t} \rVert_{F} = \sqrt{ \sum_{i, j = 1}^{2} \lVert A^{t+1}_{ij} \rVert^{2}_{F} }.
\end{equation}
Since the square of the Frobenius norm is additive in its submatrix norms we
can construct it recursively starting at the bottom tier. The matrix product is
defined recursively at each tier as
\begin{equation}
\label{eq:product}
C^{t}_{ij} = \sum_{k = 1}^{2} A^{t+1}_{ik} B^{t+1}_{kj}
\hspace{0.2in} i, j = 1, 2,
\end{equation}
subject to the \SpAMM{} condition
\begin{equation}
\label{eq:norm_condition}
A^{t+1}_{ik} B^{t+1}_{kj} \equiv 0
\hspace{0.2in} \mbox{if }
\lVert A^{t+1}_{ik} \rVert_{F} \, \lVert B^{t+1}_{kj} \rVert_{F} < \tau,
\end{equation}
where $\tau$ is the \SpAMM{} product tolerance.

\subsection{Matrix Sparsity \emph{vs.}~Product Sparsity}

Matrices with algebraic or exponential decay are not necessarily sparse unless
elements are dropped (sparsified). Storage of dense matrices is $\mathcal{O}
\left( n^{2} \right)$, while lookup and insertion/deletion is $\mathcal{O}
\left( \log n \right)$. In practice though, with $\tau > 0$, truncation in the
product space leads naturally to low level truncation in the vector space of
the result, $C$.  Likewise, it is possible to apply a very small threshold to
the product matrices $A$ and $B$ to maintain an $\mathcal{O} \left( n \log n
\right)$ computational complexity of the multiply (for further discussion, see
Sec.~\ref{sec:discussion}). In particular, using a drop tolerance $\epsilon$,
with
\begin{equation}
\epsilon = \frac{\tau}{ \max \left( \lVert A \rVert_{F}, \lVert B \rVert_{F}
\right) },
\end{equation}
is numerically consistent with the \SpAMM{} condition. The accumulation of
error due to low level sparsification and product space truncation are certain
to be application specific, with a detailed analysis beyond the scope of the
current work. Here we consider only errors associated with application of the
\SpAMM{} condition to full (dense) matrices.

\section{Implementation}
\label{sec:implementation}

\begin{figure}
\includegraphics[width=1.0\columnwidth]{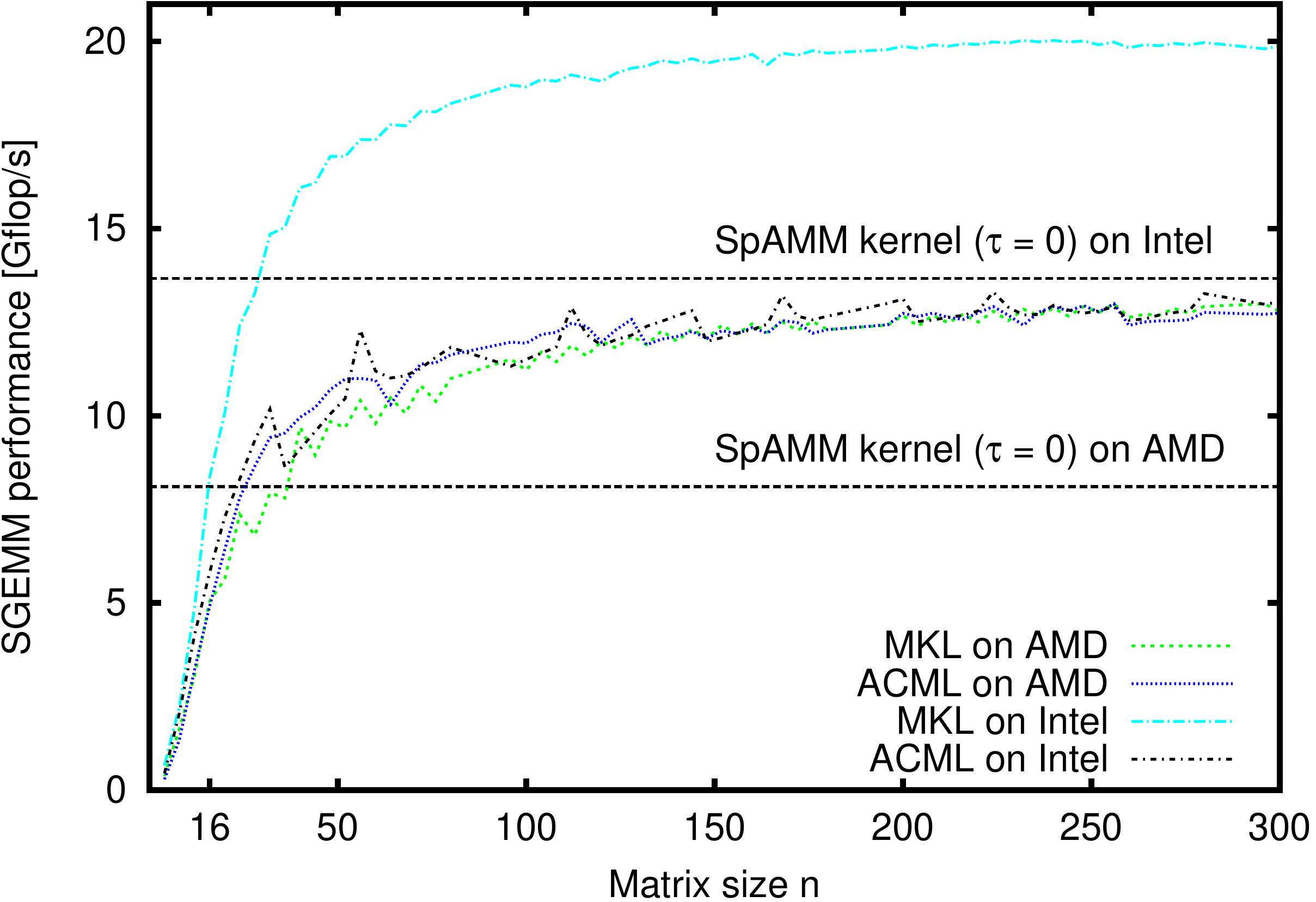}
\caption{\label{fig:sgemm_comparison} The performance of {\tt SGEMM} of the
{\tt MKL} and {\tt ACML} libraries on AMD Opteron 6168 and Intel Xeon X5650.
Clearly shown is the slow ramping up of performance with matrix size. For
comparison the average performance on a dense $16 \times 16$ submatrix of our
numeric multiply is shown as horizontal lines.}
\end{figure}

Figure \ref{fig:sgemm_comparison} shows the performance of tuned {\tt SGEMM}
kernels from {\tt MKL} \cite{intelMKL} and {\tt ACML} \cite{amdACML} with
increasing matrix size. For small matrices, the performance is mostly
throttled by memory access since memory latency is very large with $\approx$
60--100 ns (at 2.8 GHz, this translates into roughly 168--280 cycles) per
memory access to main memory \cite[Table 2]{IntelPerformance}, compared to the
latency of an addition or multiplication operation, taking only 3 or 4 cycles
\cite[Xeon X5650, Appendix C]{Intel2009}, respectively.  In order to mitigate
this latency, optimized, cache-aware {\tt SGEMM} implementations use detailed
knowledge of the cache hierarchy and data reordering techniques, such as
partitioning and packing into contiguous buffers \cite{Goto2002,
Goto:2008:AHM:1356052.1356053, Goto:2008:HIL:1377603.1377607,
gustavson2011parallel}, to replace main memory access with cache access, which
is significantly faster (access to L1 data cache takes $\approx$ 10 cycles).

\subsection{Cache-Oblivious Algorithms}

A class of alternative approaches is aimed at avoiding knowledge of the cache
hierarchy altogether through improved data locality, and is typically referred
to as cache-oblivious methods \cite{bader08hardware-oriented, Jin2005,
Wise:2001:LSM:568014.379559, Thiyagalingam03improvingthe,
Hanlon:2001:CCM:565761.565766, Adams:2006:SOS, Frigo:1999:CA:795665.796479,
OlsonSkov2002, Gottschling:2007:RMA:1274971.1274989,
Lorton:2006:ABL:1166133.1166134, Gabriel:2004:OCR:1054943.1054962,
mellor2001improving}. Typically, a heuristic such as a SFC is applied to
storage and computation ordering as shown in Fig.~\ref{fig:SFC} (C) and (D),
which yields excellent locality for both, and the newly gained locality
obviates cache-aware methods employed by standard dense methods.

Cache-oblivious algorithms effectively mitigate the issue of latency between
memory hierarchies, and would seem to be a natural way to deal with sparse
irregular problems such as \SpAMM{}. However, fine grained cache oblivious
methods neglect to address the hardware prefetcher, which can have significant
impact on application performance \cite{Pan07algorithmsto, Dick2007} and are
found to underperform cache-aware approaches for dense linear algebra problems
\cite{yotov2007experimental, gunnels2007cache}.

\subsection{Hardware Prefetch}

The prefetcher predicts the memory access pattern and pre\-loads the
anticipated addresses into cache without software intervention. Currently
though, the algorithms employed by the prefetcher are geared towards regular
contiguous problems such as dense linear algebra and require fixed stride
memory access patterns. In addition, prefetching will not cross a 4 KiB page
boundary and the stride has to be within a very short threshold value, which
depends on the processor implementation, but is typically 256 or 512 bytes.
The prefetcher is triggered after two successive cache misses in the last
level cache with address distances within the threshold value
\cite[Sec.~2.4.4.4]{Intel2009}.

In an early implementation of \SpAMM{} using cache-oblivious methods, we
concluded that the effects of the hardware prefetcher are significant.  Our
current understanding is that fine grained cache-oblivious techniques introduce
variable stride access patterns, which fail to activate the hardware prefetch
unit (two successive cache misses with constant stride are necessary), leading
to performance degradation relative to conventional cache-aware
implementations. This understanding is corroborated by the findings of Bader
\emph{et al} \cite{bader08hardware-oriented} who report poor performance of
their cache-oblivious algorithm for {\tt DGEMM} at fine granularities. While it
may be possible to improve fine-grained performance through software prefetch
statements, to the best of our knowledge this remains a challenging, unrealized
goal for irregular algorithms.

In summary, while it is desirable to apply the \SpAMM{} condition,
Eq.~(\ref{eq:norm_condition}), aggressively, fine-grained irregular memory
access precludes the use of hardware prefetch in its current form, favoring a
coarser granularity. Our kernel strikes a compromise between these
considerations by implementing storage and computation at different
granularities. We retain a quadtree structure down to a granularity of $16
\times 16$ and divide each dense block into $4 \times 4$ submatrices, stored in
row-major order but multiplied under the \SpAMM{} condition. While the required
conditionals in the kernel degrade performance only slightly in the dense case
they lead to increasingly effective performance as measured by the inverse of
time-to-solution (see Appendix \ref{sec:effective_performance}), and excellent
error control with sparsity. In the following, we detail this strategy.

\subsection{Recursive Implementation}
\label{sec:recursive_implementation}

An unoptimized version of \SpAMM{} is straightforward to write using recursion
and a conventional, vendor tuned kernel for the block multiply. For all
matrices with decay in our tests (see Sec.~\ref{sec:benchmarks}) we found the
fastest time to solution for this recursive approach and a conventional {\tt
SGEMM} kernel to occur for $16 \times 16$ blocking.

\subsection{Fast Implementation}
\label{sec:fast_implementation}

We divide the \SpAMM{} algorithm into two parts, a symbolic part that
collects the submatrices to be multiplied and a numeric part that multiplies
the submatrices, described separately in the following.

\subsubsection{Symbolic Multiply with Linkless Tree}
\label{sec:linkless_tree}

Tree data structures and recursive implementations on trees can be inefficient
in languages such as C/C++ or Fortran. While the case of tail-recursion can be
inlined by many modern compilers, the general case is significantly more
challenging and may not always allow for optimization
\cite{Stitt:2008:RF:1366110.1366143, Tang:2006:CIR:1185448.1185574}. Hashed
(linkless) tree structures offer a fast and efficient alternative. Instead of
storing explicit links, the nodes are stored in hash tables by a linearized
index, as described by Morton \cite{Morton66} and Gargantini
\cite{Gargantini1982} and the excellent book by Samet
\cite{Samet:1990:DAS:77589}. This technique is often employed in generalized
$N$-Body solvers, and in particular in the computer graphics realm for
collision detection and culling \cite{CGF:CGF1554, CGF:CGF1775,
Lefebvre:2006:PSH:1141911.1141926, Lefebvre:2006:PSH:1179352.1141926,
AVRIL:2009:HAL-00412870:1, Zou2008, 694268}, and in database applications such
as the join operation \cite{Amossen:2009:SpMMeqJoin,
Lieberman:2008:FSJ:1546682.1547260, Jacox:2003:ISJ:937598.937600, 4115860},
with some of these implemented on GPU systems.

\begin{figure}
\includegraphics[width=1.0\columnwidth]{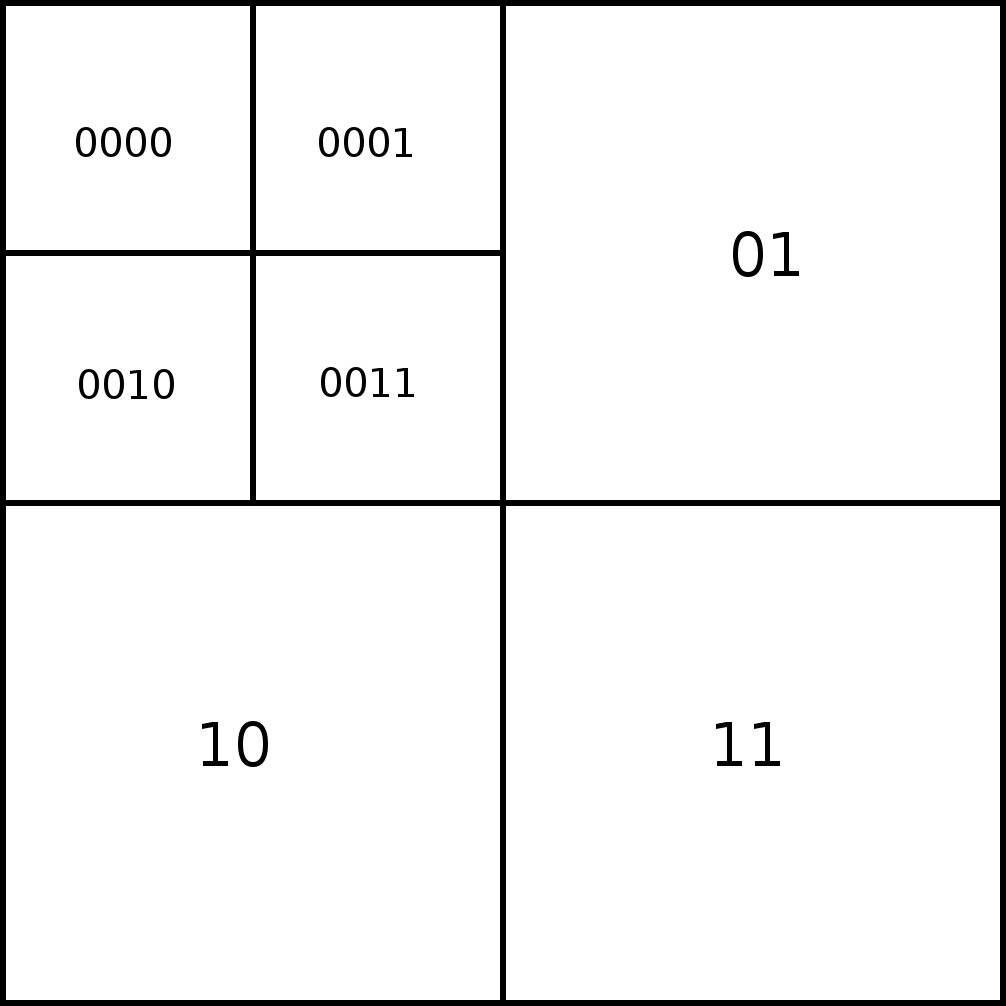}
\caption{\label{fig:linear_quadtree} A schematic representation of the
construction of linear quadtree indices and where those indices are in the
original matrix. Shown are the linear indices of first tier for the matrix
quadrants II, III, and IV. The linear indices of the second tier are shown for
quadrant I.}
\end{figure}

Each tree node is identified by its upper left hand submatrix row- and column
indices, $i$ and $j$. In combination with the matrix dimension and the tier
depth, the full range of matrix indices covered by a node can be readily
reconstructed.  The linear index $l$ is calculated through dilating and
interleaving $i$ and $j$:
\begin{eqnarray}
i     & = & \sum_{m = 0} i_{m} 2^{m} \\
j     & = & \sum_{m = 0} j_{m} 2^{m} \\
\label{eq:linear_index}
l     & = & \sum_{m = 0} \left[ i_{m} 2^{2m+1} + j_{m} 2^{2m} \right].
\end{eqnarray}
For illustration, Fig.~\ref{fig:linear_quadtree} shows the first two tiers of
linear matrix indices of a small matrix. Notice how each tier adds two digits
from the right to the linear index of its parent node and the quadrant labeling
scheme is recursively repeated. Despite the lack of links between tree nodes,
the full tree structure and its hierarchy can simply be reconstructed on the
way down by adding two digits from the right, or on the way up by removing
those two digits.

It is worth pointing out that for performance the choice of hash function is
extremely important. Excessive collision rates are to be avoided, since hash
table performance degrades rapidly with increasing collision rate. We use the
well designed 32-bit hash function of Jenkins \cite{Jenkins2006}, which limits
the tree depth to $16+4$ (since we store $16 \times 16$ dense submatrices) and
the matrix size to $1048576 \times 1048576$, more than sufficient for the tests
presented here. Larger matrices are possible through the use of either a larger
hash key or a hybrid approach of linked tree (at the top) and linkless tree (at
the bottom 20 tiers).

\subsubsection{Symbolic Multiply}
\label{sec:symbolic_multiply}

The convolution of the two-dimensional vector spaces into the
three-dimensional product space, panels (C) and (D) of Fig.~\ref{fig:SFC}, is
straightforward to implement through recursively applying
Eq.~(\ref{eq:product}), but the situation is more complicated in the case of
linkless trees. The hash tables of $A$ and $B$ have to be searched and
matching submatrix pairs $A_{ik}$ and $B_{kj}$ have to be found, conditional
on Eq.~(\ref{eq:norm_condition}), which result in a contribution to the
resulting submatrix $C_{ij}$. Related problems are common amongst generalized
$N$-Body solvers, such as performing queries on hashed trees to join
\cite{Amossen:2009:SpMMeqJoin}, collide \cite{AVRIL:2009:HAL-00412870:1}, add
\cite{Warren:1995:HOT}, or in other ways combine elements.

The collection of the submatrices that contribute to the product is carried
out in several steps. First, the linear indices $l$ for $A$ and $B$ are
extracted from the hash tables of the lowest tier ($16 \times 16$ submatrix
level) and stored in an array. From Eq.~(\ref{eq:linear_index}) the indices
are given bit-wise by:
\begin{equation}
l_{A, ik} = \ldots i_{2} k_{2} i_{1} k_{1} i_{0} k_{0}.
\end{equation}

\noindent Second, in a first pass, the index arrays are sorted on their
$k$-index values using the merge sort algorithm\footnote{Although quicksort is
often found to be faster than merge sort in practice, it turns out to be a
poor choice in our case: Sorting $n$ elements, merge sort exhibits average and
worst-case cost of $\mathcal{O} \left(n \log n \right)$, while in comparison
quicksort exhibits the same average cost, but slows down to $\mathcal{O}
\left( n^{2} \right)$ in the worst-case. For the matrices we used in our
tests, quicksort tended to deviate from logarithmic towards quadratic
behavior, and merge sort consistently performed better.  We speculate that a
different choice of hash function could improve the performance of quicksort
and the sorting steps in the symbolic part, but the design of a good hash
function which exhibits few collisions and preserves logarithmic quicksort
behavior is beyond the scope of this work.}, resulting in an ordered array
partitioned into blocks with identical $k$-values; ``$k$-value blocks''.
Application of the bit-mask {\tt 0b\dots010101} to $l_{A, ik}$ and {\tt
0b\dots 101010} to $l_{B, kj}$ renders this step straightforward. In a second
pass, within each $k$-value block, the indices are sorted by their associated
node norms (the Frobenius norm stored at each node, see Sec.~\ref{sec:spamm}
for details) in descending order.  Alternatively both sort passes can be
combined into one pass by rewriting the sort comparison function to include
both $k$-index and norm, however we have not found any performance benefit
from this approach.

\begin{algorithm}
\caption{\label{algorithm:spamm_multiply} Linkless tree \SpAMM{} multiply.}
\algsetup{indent=2em}
\begin{algorithmic}[1]
\STATE Extract all keys (linear matrix indices) from $A$ and $B$.
\STATE Sort key lists on $k$, {\it i.e.} $A_{ik}$ and $B_{kj}$.
\FORALL{$k$ in key lists}
  \STATE Sort sublist by descending norm.
\ENDFOR
\STATE Convolve and find $A$, $B$, and $C$ submatrices,
Eq.~(\ref{eq:convolution}).
\FORALL{products $\in$ product array}
  \STATE $C_{ij} \leftarrow C_{ij} + A_{ik} B_{kj}$
\ENDFOR
\end{algorithmic}
\end{algorithm}

\noindent Third, the convolution is performed by nested loops over the
$k$-value blocks of the arrays. Although our algorithm only considers the
lowest tier, rendering it non-hierarchical, norm sorting the $k$-value blocks
restores $\mathcal{O} \left( n \log n \right)$ or $\mathcal{O} \left( n^{2}
\log n^{2} \right)$ complexity for sparse and dense matrices, respectively.
We implemented two versions: The first consists of straightforward nested
loops (``nested convolution''), while the second is partially unrolled and
uses Intel's SSE intrinsics (``unrolled convolution'') that leads to a
speed-up of close to three. Since SSE instructions require proper alignment of
memory, we align each $k$-value block on a 16-byte boundary and the index
computations, the norm products, and the product comparisons are done in SSE,
which reduces the loops to a stride of four. We note that the index extraction
from the hash tables leads to $\mathcal{O} \left( n^{2} \right)$ storage for
the indices of dense matrices.  While this is not prohibitive for the matrices
tested here, it may become so for larger matrices if not ultimately
sparsified.  Pointers to the nodes of each triple of submatrices $A_{ik}$,
$B_{kj}$, and $C_{ij}$ are then extracted from the hash tables and appended to
an array for further processing by our numeric multiply, described in the
following section. The $C_{ij}$ submatrix indices are constructed by masking
out the $i$ and $j$ index values of $l_{A, ik}$ and $l_{B, kj}$ and combining
the two masked indices,
\begin{equation}
\label{eq:convolution}
l_{C, ij}  =
\left( l_{A, ik} \,\, \& \,\, \mbox{0b\dots101010} \right)
|
\left( l_{B, kj} \,\, \& \,\, \mbox{0b\dots010101} \right).
\end{equation}
Steps 1--6 of Algorithm \ref{algorithm:spamm_multiply} show pseudocode of the
convolution step.

While storing the nodes of the linkless tree in a hashtable is very efficient
for sparse matrices, a simple array is faster for dense matrices. Since the
matrices used in this work are all dense (albeit with the aforementioned decay
properties), it is educational to consider a performance comparison between
hashtable and array storage of the indices. We therefore benchmarked two
implementations of the linkless tree: Hashtable storage of matrix nodes with
the nested convolution implementation (labeled ``hashtable'' in
Fig.~\ref{fig:symbolic_cost}) and array storage of matrix nodes with the
unrolled convolution implementation using SSE intrinsics (labeled ``SSE'' in
Fig.~\ref{fig:symbolic_cost}).

\subsubsection{Numeric Multiply}
\label{sec:fast_kernel}

Steps 7--9 of Algorithm \ref{algorithm:spamm_multiply} show pseudocode for the
numeric multiply, and how it serially processes submatrix products.  Several
factors need to be considered for optimal performance of the kernel: On the
one hand, Streaming SIMD Extension (SSE) instructions exhibit higher
performance than comparable x87 instructions \cite[Sec.~3.8.4]{Intel2009}) and
in single-precision the SSE vector length is 4, suggesting a granularity of $4
\times 4$ for applying Eq.~(\ref{eq:norm_condition}).  On the other hand, this
$4 \times 4$ level of granularity leads to inefficient cache use and memory
access which can only be remedied by a larger contiguous data structure.
Through experimentation, we found $16 \times 16$ contiguous blocks to give the
best overall performance. Finally, optimal data storage and access patterns
had to be identified. We considered row-major, Morton order, and hierarchical
versions of these and found that simple non-hierarchical row-major storage and
access is the most efficient.

A note on the hardware prefetch: The $16 \times 16$ single-precision
submatrices occupy 16 cache lines (or more as in our SSE version described
below), which is sufficiently large for the hardware prefetch to accelerate
memory access. The conditionals around the $4 \times 4$ block multiplies do not
significantly degrade performance by themselves, but when products are skipped
due to Eq.~(\ref{eq:norm_condition}), memory access becomes irregular and
prevents the hardware prefetcher from activating, which degrades the flop rate.
At the same time however, the effective performance, the inverse of
time-to-solution, increases since more products are dropped, more than
offsetting the loss of hardware prefetch.

\begin{algorithm}
\caption{\label{algorithm:SSE_version} SSE version of $4 \times 4$ multiply.}
\begin{algorithmic}[1]
\STATE register[1] $\leftarrow [ 0 \, 0 \, 0 \, 0 ]$
\STATE register[2] $\leftarrow [ A_{11} \, A_{11} \, A_{11} \, A_{11} ]$
\STATE register[3] $\leftarrow [ B_{11} \, B_{12} \, B_{13} \, B_{14} ]$
\STATE register[2] $\leftarrow$ register[2] $\circ$ register[3]
\STATE register[1] $\leftarrow$ register[1] $+$ register[2]
\STATE register[2] $\leftarrow [ A_{12} \, A_{12} \, A_{12} \, A_{12} ]$
\STATE register[3] $\leftarrow [ B_{21} \, B_{22} \, B_{23} \, B_{24} ]$
\STATE register[2] $\leftarrow$ register[2] $\circ$ register[3]
\STATE register[1] $\leftarrow$ register[1] $+$ register[2]
\STATE register[1] $\rightarrow [ C_{11} \, C_{12} \, C_{13} \, C_{14} ]$
\end{algorithmic}
\end{algorithm}

We implemented two versions of the kernel in assembly using the SSE and the
SSE4.1 instruction sets. The Xeon X5650 supports both instruction sets, while
the Opteron 6168 only has support for SSE.  The SSE instruction set was
introduced by Intel in 1999 with their Pentium III processor series, adding 70
instructions and eight registers (sixteen on 64-bit processors).  Lacking a
dot product, our SSE kernel relies on a dilated storage scheme for the matrix
elements of $A$, in which each matrix element is stored 4 times, {\it i.e.} in
vectors such as $[ A_{11} \, A_{11} \, A_{11} \, A_{11} ]$. This storage
scheme allows for the efficient use of SSE registers without having to shuffle
or copy the result matrix elements in SSE, leading to a performance
improvement despite the increased demand on memory transfer.  We show
pseudocode of this kernel version in Algorithm \ref{algorithm:SSE_version}.
Note that the SSE multiplication and addition instructions used here operate
element-wise with the Hadamard or Schur product denoted by the symbol $\circ$.

The SSE4 instruction set was introduced by Intel with their Penryn processor
core in 2006. A subset of SSE4, SSE4.1, introduces dot product instructions for
packed single-precision numbers and the dilated storage scheme of the SSE
version becomes unnecessary. The matrix elements of $B$ are stored in
column-major order.  Compared to the SSE version, SSE4.1 allows for a
reduction of data movement and it uses a dedicated dot product instruction.
However, we did not find this version to be faster than the SSE version on the
Intel Xeon system.

\section{Experimental Methodology}
\label{sec:methodology}

In the following, we give a detailed description of our experimental
methodology. In Sec.~\ref{sec:platforms} we will describe in detail the
hardware platforms used, and in Sec.~\ref{sec:benchmarks} the benchmarks used
to assess the performance of our implementation.

\subsection{Platforms}
\label{sec:platforms}

All performance measurements were conducted on two hardware platforms, Intel
Xeon X5650 (Westmere microarchitecture, 2.66 GHz, 6 cores, 32 KiB of L1 data
per core, 256 KiB of L2 per core, 12 MiB L3 shared by all cores) and AMD
Opteron 6168 (Magny-Cours microarchitecture, 1.8 GHz, 12 cores, 64 KiB L1 data
per core, 512 KiB L2 per core, 12 MiB L3 shared by all cores), using a NUMA
aware Linux kernel ({\tt CONFIG\_NUMA} and {\tt CONFIG\_X86\_64\_ACPI\_NUMA})
version 2.6.39 on Gentoo \cite{gentoo}.  Cycle and flop counts where taken
with the {\tt PAPI} library \cite{papi, Mucci99papi:a, 5936277} version
4.1.4\footnote{Cycle count: {\tt PAPI\_TOT\_CYC}; Flop count, AMD: {\tt
RETIRED\_SSE\_OPERATIONS:OP\_TYPE:ALL}, Intel: {\tt
FP\_COMP\_OPS\_EXE:SSE\_FP\_PACKED}, {\tt FP\_COMP\_OPS\_EXE:SSE\_FP\_SCALAR};
Cache misses: {\tt DATA\_CACHE\_MISSES}}. Throughout, gcc \cite{gcc} version
4.4.5 was used for compilation. It should be noted that the choice of compiler
has very little impact on the results of this study because all time intensive
computational kernels either are linked in from externally compiled libraries
({\tt MKL} and {\tt ACML}) or are written in assembly.  Care has to be taken
for getting accurate and reproducible measurement results \cite{Weaver2008};
the benchmark process was locked to a processor core and memory allocation was
limited to the local NUMA node using the {\tt numactl} command in combination
with the {\tt --physcpubind} and {\tt --membind} command line options. In
addition all allocated memory was pinned to prevent paging using the {\tt
mlockall()} system call.

\subsection{Benchmarks}
\label{sec:benchmarks}

In a typical quantum chemistry calculation, the solution of an
integro-differential equation is sought by expanding the electronic
wavefunctions of a molecule such as the one shown in panel (A) of
Fig.~\ref{fig:SFC} into a set of atomic basis functions yielding a matrix
eigenvalue problem \cite{szabo1996modern}, or spectral projection
\cite{McWeeny12061956, PhysRevB.58.12704, Challacombe:1999:DMM,
Challacombe:2000:SpMM, Benzi:2007:Decay, BenziDecay2012}.  The exact form and
number of functions in the basis set and hence the number of expansion
coefficients per atom depend on the atom type and the basis set.  Because SFC
ordering is applied atom-wise, matrix elements are clustered by magnitude,
grouped in submatrices given by the atom type and basis set. We explore the
effect of submatrix size with different basis sets, STO-2G and 6-31G${}^{**}$,
for quantum chemical matrices corresponding to Hilbert curve ordered water
clusters. These basis sets introduce submatrices of different sizes (STO-2G: H
$1 \times 1$, O $5 \times 5$; 6-31G${}^{**}$: H $5 \times 5$, O $15 \times
15$), representing the native granularities embedded in the matrix structures.
The matrices were generated at the Restricted Hartree-Fock (RHF) level of
theory with {\tt FreeON}, a suite of programs for $\mathcal{O} (N)$ quantum
chemistry \cite{FreeON} and are fully converged and without sparsification
($\epsilon = 0$), {\it i.e.} the matrices are fully dense but exhibit the
approximately exponential decay pattern shown in
Fig.~\ref{fig:sparsity_pattern}. The sequence of water clusters corresponds to
standard temperature and pressure and have been used in a number of previous
studies \cite{Challacombe:Review, Challacombe:1996:QCTCb,
Challacombe:1997:QCTC, burant:8969, ESchwegler97, millam:5569, daniels:425,
ochsenfeld:1663}.

The test calculations consist of taking the matrix square, a key step in
construction of the spectral projector (density matrix purification)
\cite{McWeeny12061956, PhysRevB.58.12704, Challacombe:1999:DMM,
Niklasson:2003:TRS4, Niklasson:2004:DMPT}.  The error of the product is
measured by the max norm of the difference to the double-precision ({\tt
DGEMM}) product:
\begin{equation}
\label{eq:matrix_error}
\left\| \Delta C \right\|_{\mathrm{max}} = \max_{ij} \left\{ \left| C_{ij} -
C^{\mbox{\tt DGEMM}}_{ij} \right| \right\},
\end{equation}
which places an upper bound on the element-wise error and is a simple check on
the accuracy of the multiplication. For the performance comparisons the
tolerance was adjusted to match the errors. While the exact dependence of the
product error on $\tau$ is currently unknown, it is worth pointing out that
due to its application in product space, $\tau$ affects which product
contributions are dropped linearly, as opposed to the quadratic dependence on
$\epsilon$ for sparsification \cite{ChallacombeBock2010}. Performance is
reported using three metrics, the measured flop rate, the effective
performance, and the cycle count as reported by {\tt PAPI}. For reference, the
effective performance is proportional to the inverse of the execution time,
see Appendix \ref{sec:effective_performance} for details.

\begin{table}
\caption{\label{table:benchmarks} Performance and error benchmarks}
\begin{center}
\begin{tabular}{rp{0.55\columnwidth}ll}
\hline
\hline
\multicolumn{2}{c}{Test} & \hspace{0.2cm} & Label \\
\hline
1.
  & Comparison of convolution implementations
  &
  & Symbolic Multiply \\
2.
  & Dense multiply with {\tt SGEMM} from {\tt MKL} and {\tt ACML}
  (Fig.~\ref{fig:sgemm_comparison})
  &
  & {\tt SGEMM} \\
3.
  & \SpAMM{} numeric multiply at $4 \times 4$
  &
  & \SpAMM{}($4 \times 4$) \\
4.
  & \SpAMM{} numeric multiply at $16 \times 16$, no conditionals
  &
  & \SpAMM{}($16 \times 16$) \\
5.
  & \SpAMM{} multiply at $16 \times 16$ level, {\tt SGEMM}
  &
  & \SpAMM{}({\tt SGEMM}) \\
\hline
\hline
\end{tabular}
\end{center}
\end{table}

Test 1 consists of a comparison between the different convolution
implementations, described in Sec.~\ref{sec:symbolic_multiply}.  Test 2
consists of matrix multiplications with vendor tuned {\tt SGEMM} for varying
matrix sizes to asses the performance of our numeric multiply
(Sec.~\ref{sec:fast_kernel}).  The other tests consist of matrix
multiplications with \SpAMM{} to assess the performance and product error of
\SpAMM{}. All \SpAMM{} tests use the unrolled convolution implementation with
SSE intrinsics, as described in Sec.~\ref{sec:symbolic_multiply}, for the
symbolic multiply since it proved to be the fastest, and different versions of
the numeric multiply. The numeric multiply used in test 3 applies the \SpAMM{}
condition at a $4 \times 4$ submatrix level, which we consider our target
granularity. The kernel in test 4 uses a version of the numeric multiply in
which the \SpAMM{} conditionals are commented out, which effectively means
that the \SpAMM{} condition is applied at a $16 \times 16$ submatrix level.
Note that the execution flow is identical to test 3. The numeric multiply used
in test 5 is a vendor tuned {\tt SGEMM} at a $16 \times 16$ submatrix level. A
summary of all benchmark tests can be found in Table \ref{table:benchmarks}.

\section{Results}
\label{sec:results}

Our results are qualitatively identical on the AMD and the Intel platforms and
we will present only the AMD results in the following. In addition, we did not
find any difference in the error between {\tt MKL} and {\tt ACML} and report
for the error, Figs.~\ref{fig:water_631gss_RHF_STO_2G_good_error_AMD} and
\ref{fig:water_631gss_RHF_6_31Gss_tight_error_AMD}, only the {\tt MKL}
results.

\subsection{Symbolic Multiply}

\begin{figure}
\includegraphics[width=1.0\columnwidth]{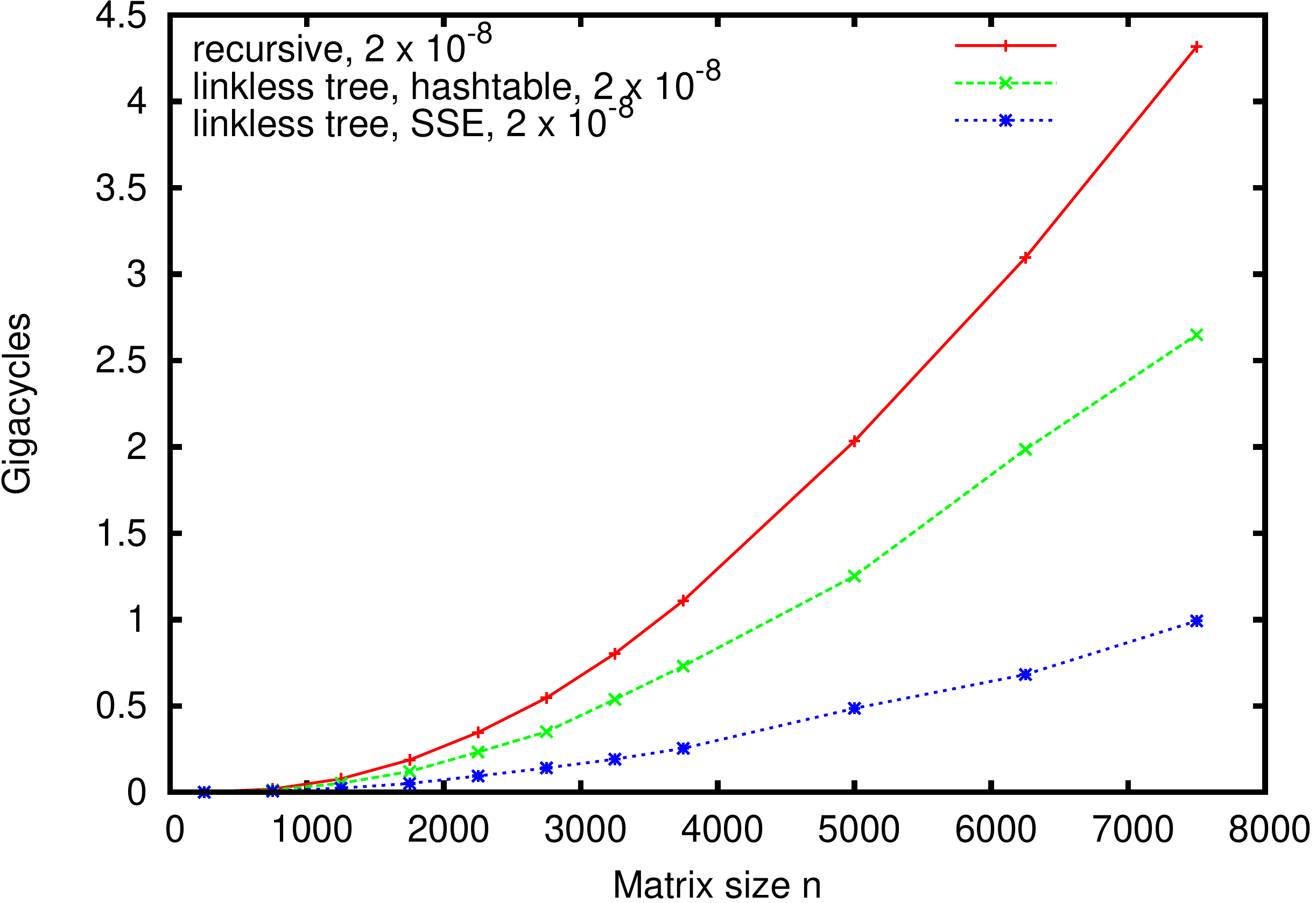}
\caption{\label{fig:symbolic_cost} The cycle counts of the symbolic part of
\SpAMM{} with a tolerance of $\tau = 2 \times 10^{-8}$ for RHF/6-31G${}^{**}$
water cluster for the recursive multiply and the two linkless tree
implementations.}
\end{figure}

Fig.~\ref{fig:symbolic_cost} (test 1 of Table \ref{table:benchmarks}) shows a
cycle count comparison of the symbolic part of \SpAMM{}($4 \times 4$, $\tau =
2 \times 10^{-8}$) between the recursive and the two linkless tree versions
(see Sec.~\ref{sec:symbolic_multiply}). For a matrix of $n = 7500$, the
unrolled convolution using arrays is almost 4.5 times faster than the
recursive version and about 2.5 times faster than the nested convolution using
hashtables. At this matrix size, the symbolic part accounts for about 12\% of
the total cost of the multiply.

\subsection{Numeric Multiply \emph{vs.}~{\tt SGEMM}}

Figure \ref{fig:sgemm_comparison} (test 2 of Table \ref{table:benchmarks})
shows the performance of {\tt SGEMM} on the Intel and AMD platforms from {\tt
MKL} and {\tt ACML}. The performance slowly ramps up and asymptotes to 12.6
Gflop/s for {\tt ACML} and {\tt MKL} on AMD and 12.7 Gflop/s for {\tt ACML}
and 19.8 Gflop/s for {\tt MKL} on Intel.  Also shown in the figure the average
performance of our numeric multiply on a dense $16 \times 16$ block for large
matrices, exhibiting a performance of 8.1 (64\% of peak {\tt ACML}/{\tt MKL})
and 13.7 (69\% of peak {\tt MKL}) Gflop/s on AMD and Intel, respectively. At a
matrix size of $16 \times 16$, the performance ratio between the \SpAMM{}
kernel and {\tt SGEMM} is 162\% (5.0 Gflop/s for {\tt MKL} and {\tt ACML}) for
AMD and 161\% (8.5 for {\tt MKL}) for Intel.

\subsection{Errors}

\begin{figure}
\includegraphics[width=1.0\columnwidth]{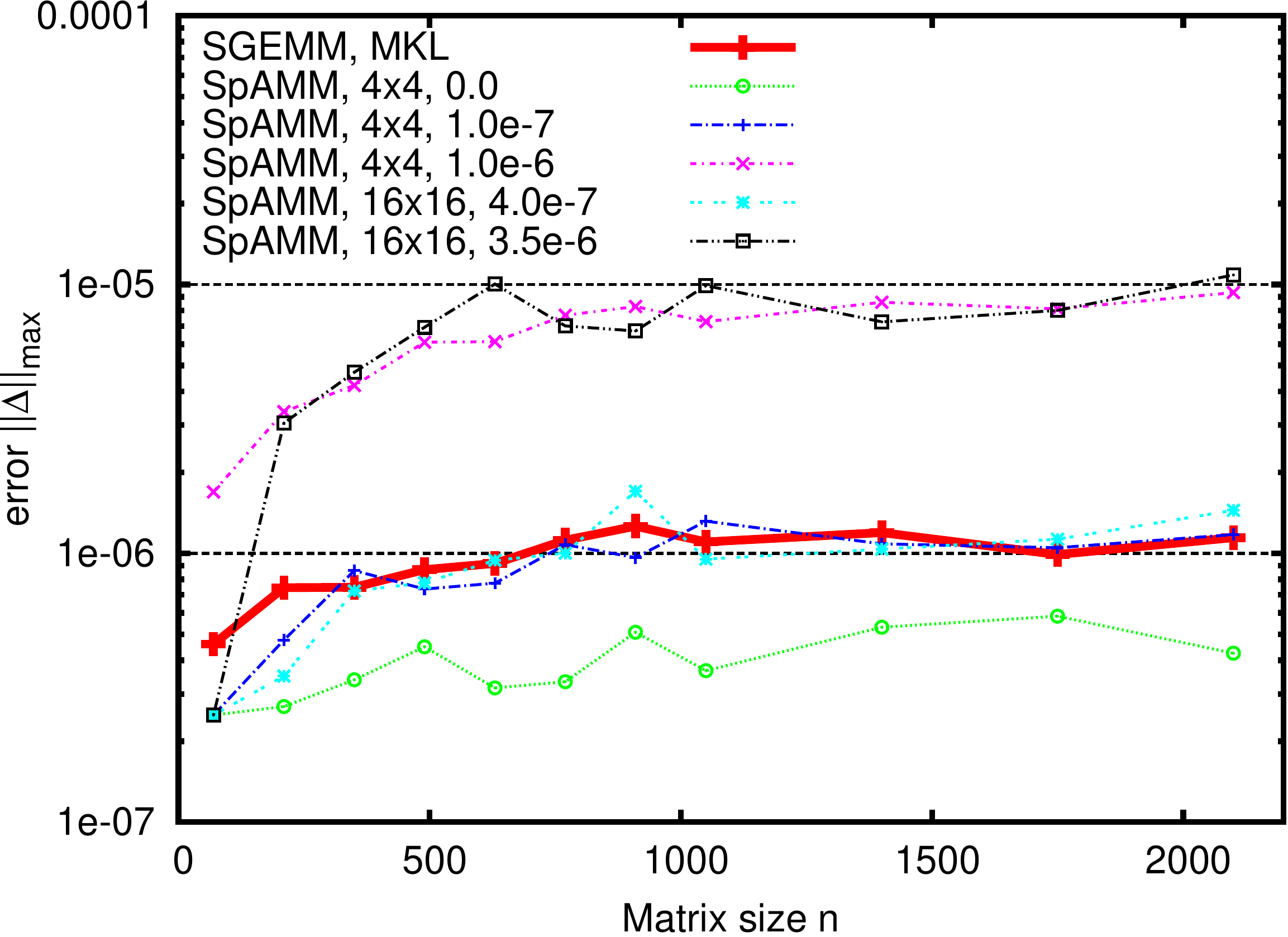}
\caption{\label{fig:water_631gss_RHF_STO_2G_good_error_AMD} Error comparison as
measured by the max norm, Eq.~(\ref{eq:matrix_error}), for water clusters of
different sizes in RHF/STO-2G on AMD.}
\end{figure}

\begin{figure}
\includegraphics[width=1.0\columnwidth]{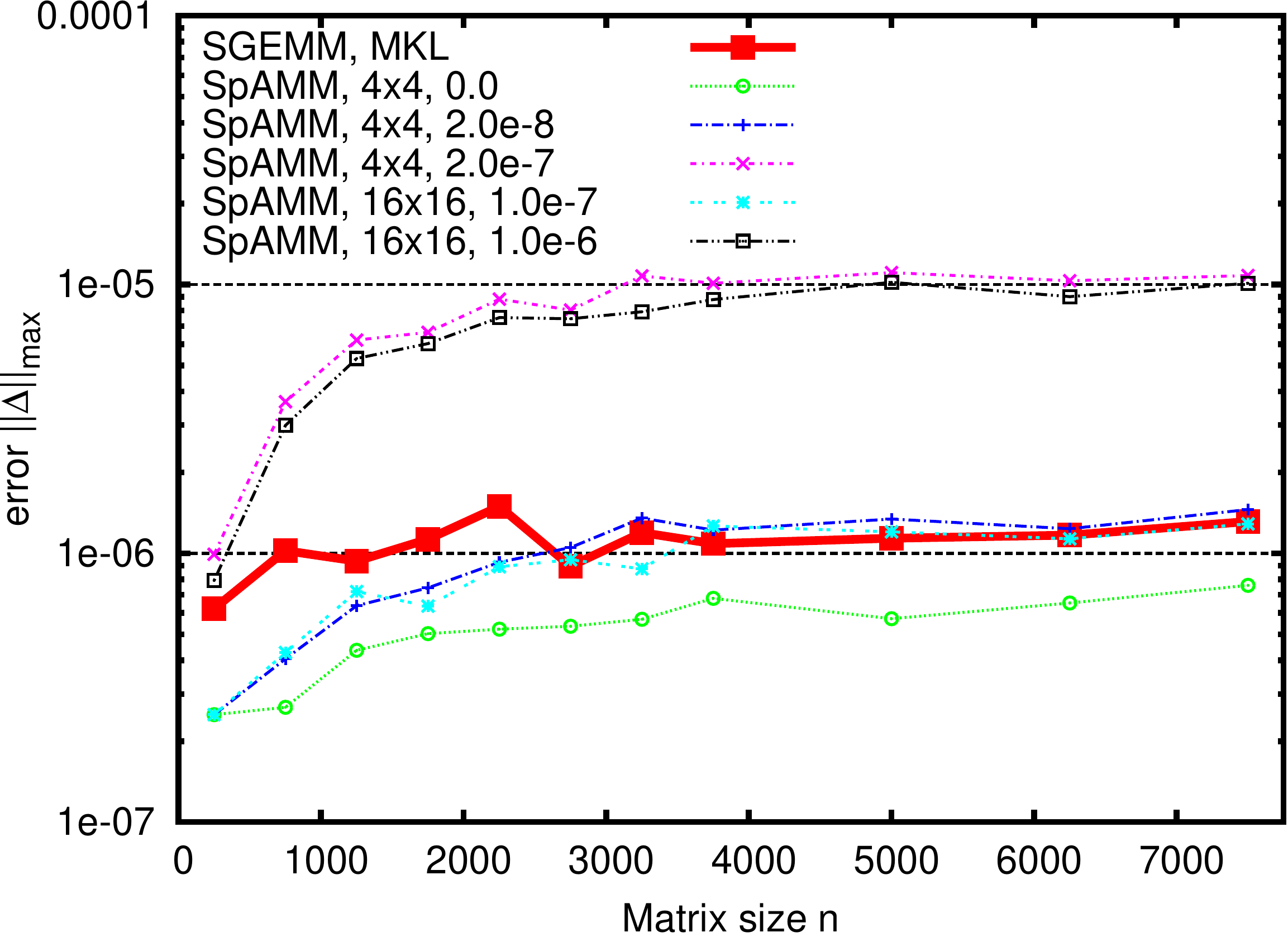}
\caption{\label{fig:water_631gss_RHF_6_31Gss_tight_error_AMD} Error comparison
as measured by the max norm, Eq.~(\ref{eq:matrix_error}), for water clusters of
different sizes in RHF/6-31G${}^{**}$ on AMD.}
\end{figure}

In order to facilitate a fair comparison between tests 3, 4, and 5
(Table \ref{table:benchmarks}), appropriate tolerance values were found to
match the matrix errors, Eq.~(\ref{eq:matrix_error}), in the different tests.
In Figs.~\ref{fig:water_631gss_RHF_STO_2G_good_error_AMD} and
\ref{fig:water_631gss_RHF_6_31Gss_tight_error_AMD} the matrix errors of the
test calculations for the two basis sets are shown. The differences in error
between tests 4 and 5 were found to be insignificant and results of test 5 are
not shown for clarity in these figures.  For comparison, a dense multiply
exhibits an error of around $10^{-6}$ and we note that \SpAMM{}($4 \times 4$,
$\tau = 0$) achieves a smaller error than {\tt SGEMM} for both basis sets,
which is likely due to differences in execution order; hierarchical
\emph{vs.}~row-column.  For the following tests of \SpAMM{} we chose two error
targets, $10^{-6}$ and $10^{-5}$, and found the tolerance values shown in
Table \ref{table:tolerance_values}. Increasing the granularity from $4 \times
4$ to $16 \times 16$ requires an increase of the tolerance by about a factor
of 4 for the STO-2G basis set, and a factor of 5 for the 6-31G${}^{**}$ basis
set.

\begin{table}
\caption{\label{table:tolerance_values} The tolerance values identified for the
water clusters at different \SpAMM{} granularity and errors.}
\begin{center}
\begin{tabular}{lllcl}
\hline
\hline
basis set                       & \hspace{0.2cm} & error                      & granularity    & tolerance \\
\hline
\multirow{4}{*}{STO-2G}         &                & \multirow{2}{*}{$10^{-6}$} & $4  \times  4$ & $1.0 \times 10^{-7}$ \\
                                &                &                            & $16 \times 16$ & $4.0 \times 10^{-7}$ \\
\cline{3-5}
                                &                & \multirow{2}{*}{$10^{-5}$} & $4  \times  4$ & $1.0 \times 10^{-6}$ \\
                                &                &                            & $16 \times 16$ & $3.5 \times 10^{-6}$ \\
\cline{1-5}
\multirow{4}{*}{6-31G${}^{**}$} &                & \multirow{2}{*}{$10^{-6}$} & $4  \times  4$ & $2.0 \times 10^{-8}$ \\
                                &                &                            & $16 \times 16$ & $1.0 \times 10^{-7}$ \\
\cline{3-5}
                                &                & \multirow{2}{*}{$10^{-5}$} & $4  \times  4$ & $2.0 \times 10^{-7}$ \\
                                &                &                            & $16 \times 16$ & $1.0 \times 10^{-6}$ \\
\hline
\hline
\end{tabular}
\end{center}
\end{table}

\subsection{Effective Performance of {\tt SGEMM} \emph{vs.}~\SpAMM{}}

\begin{figure}
\includegraphics[width=1.0\columnwidth]{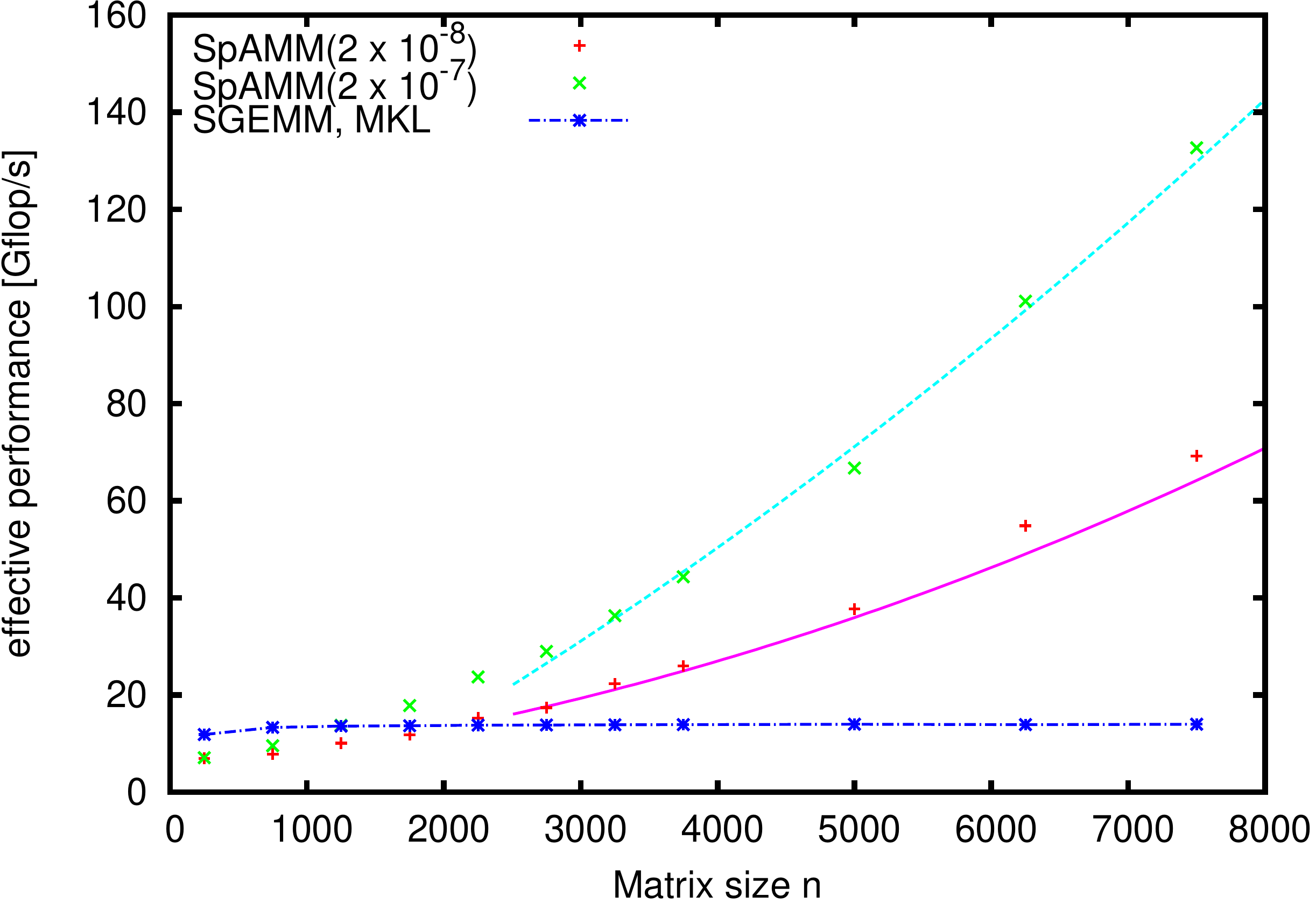}
\caption{\label{fig:water_631gss_sgemm_comparison_AMD} Effective performance of
\SpAMM{}($4 \times 4, \tau = 2 \times 10^{-8}$) and \SpAMM{}($4 \times 4, \tau
= 2 \times 10^{-7}$) and {\tt SGEMM} for the water clusters in
RHF/6-31G${}^{**}$ on AMD. Solid lines show fits to second order polynomial.}
\end{figure}

Comparing the effective performance of {\tt SGEMM}, \SpAMM{}($4 \times 4$,
$\tau = 2 \times 10^{-8}$) and \SpAMM{}($4 \times 4$, $\tau = 2 \times
10^{-7}$), for the RHF/6-31G${}^{**}$ water clusters in
Fig.~\ref{fig:water_631gss_sgemm_comparison_AMD}, we find an approximately
quadratic performance increase of \SpAMM{} with respect to matrix size. Fits to
a quadratic are shown as solid lines. The {\tt SGEMM} performance is matrix
size independent above $n \approx 1000$ due to memory bandwidth saturation and
the $\mathcal{O} \left( n^{3} \right)$ complexity of the standard matrix
multiply. We note the early cross-over of \SpAMM{} at around $n \sim 2000$ and
the fact that \SpAMM{} performs at approximately 50\% of peak {\tt SGEMM}
performance as $n \rightarrow 0$ consistent with the results shown in
Fig.~\ref{fig:sgemm_comparison}.

\subsection{Performance of \SpAMM{}}

\begin{figure}
\includegraphics[width=1.0\columnwidth]{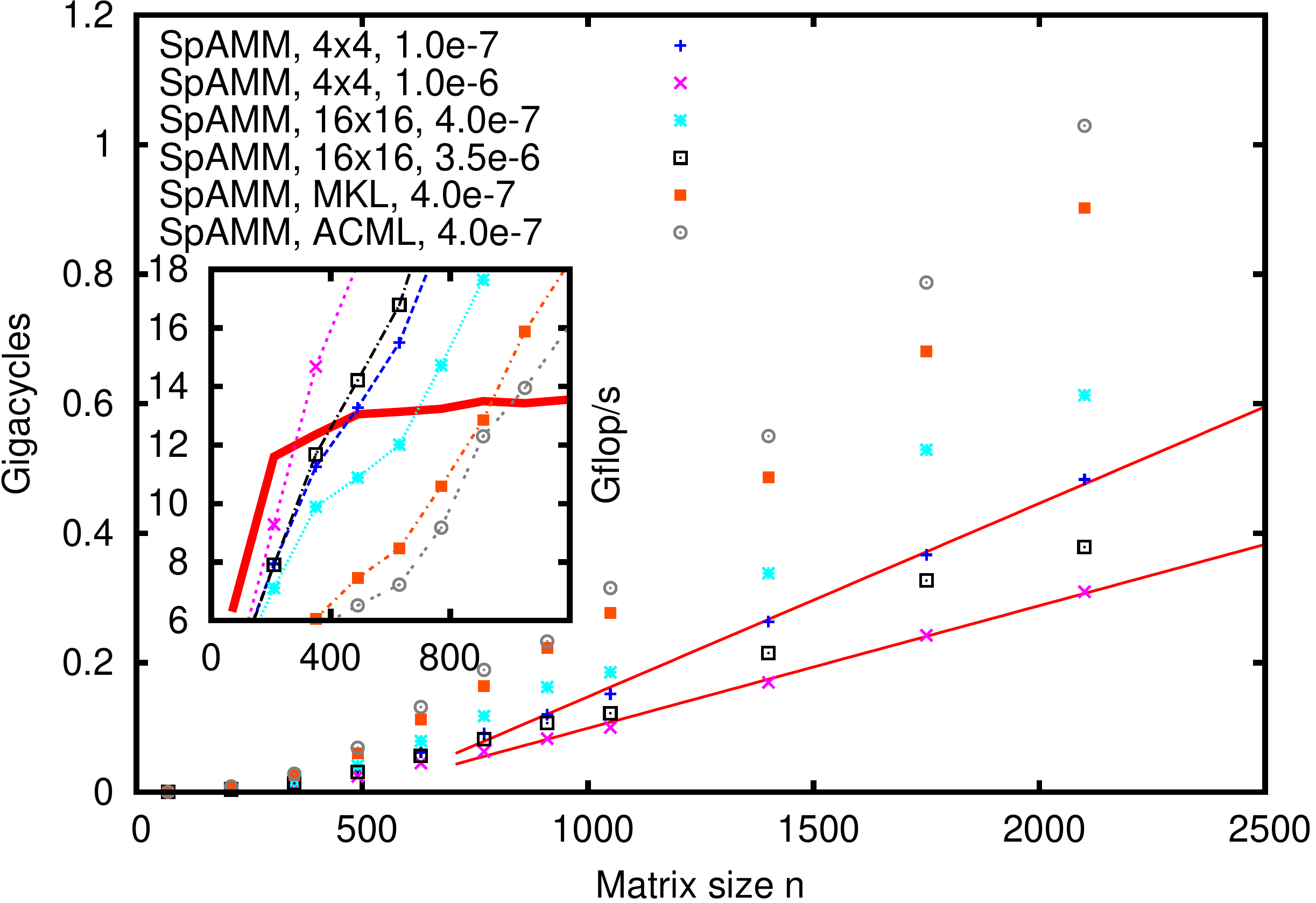}
\caption{\label{fig:water_STO2G_AMD} Performance comparison for water clusters
of different sizes in RHF/STO-2G on AMD. Solid lines show fits to a line.
Shown in the inset, the effective performance in Gflop/s compared to {\tt
SGEMM} from {\tt MKL} and {\tt ACML} (thick line).}
\end{figure}

\begin{figure}
\includegraphics[width=1.0\columnwidth]{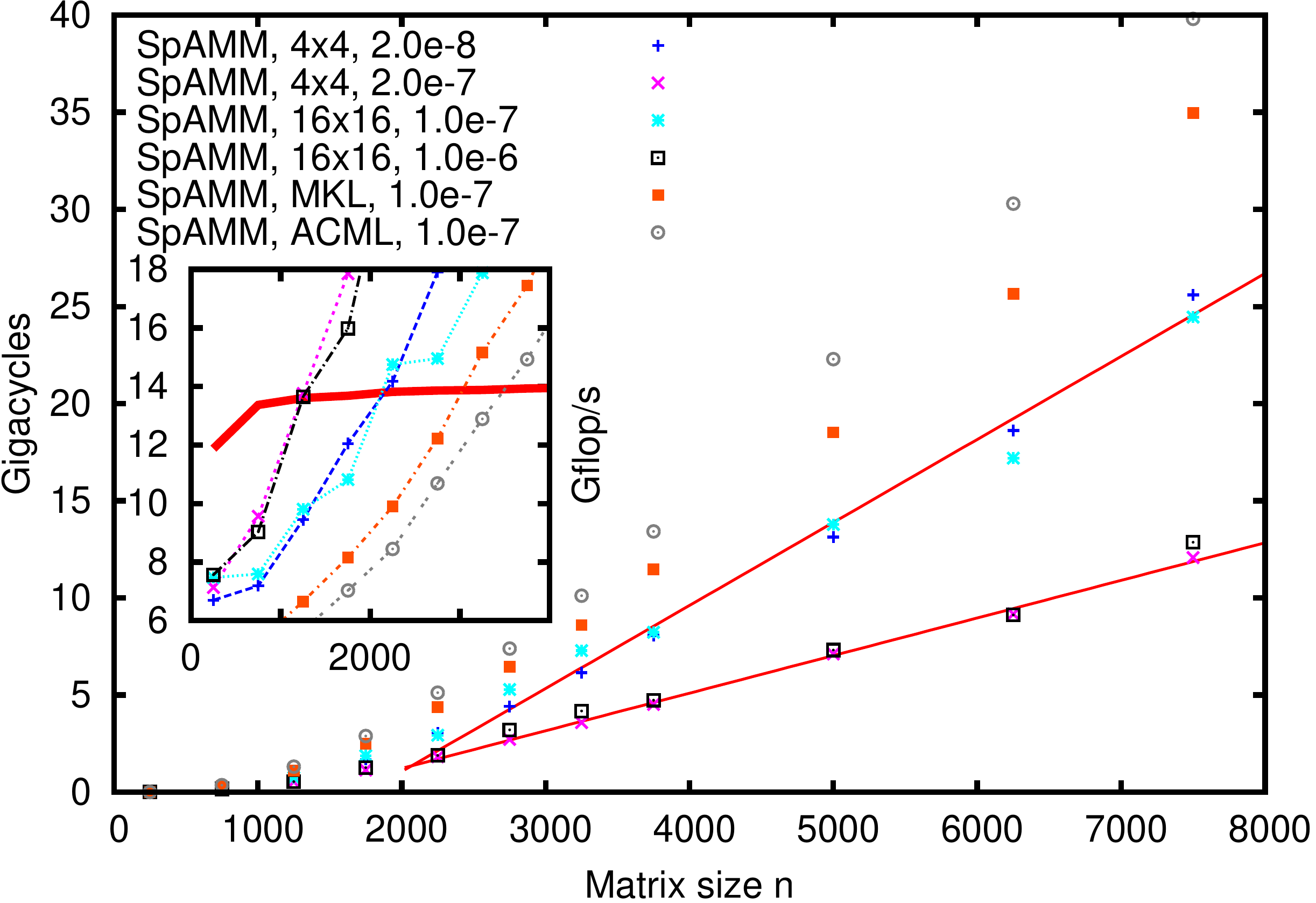}
\caption{\label{fig:water_631gss_AMD} Performance comparison for water clusters
of different sizes in RHF/6-31G** on AMD. Solid lines show fits to a line.
Shown in the inset, the effective performance in Gflop/s compared to {\tt
SGEMM} from {\tt MKL} and {\tt ACML} (thick line).}
\end{figure}

Figures \ref{fig:water_STO2G_AMD} and \ref{fig:water_631gss_AMD} show the
measured performance of tests 3, 4, and 5 of Table \ref{table:benchmarks} for
the tolerance values of Table \ref{table:tolerance_values}. Fits to a line are
shown to highlight the approximate linear scaling behavior of \SpAMM{} and to
give an indication of slope. The cycle count, which is proportional to CPU
time, is shown in the main graph, and the effective performance in the inset.
For comparison the performance of {\tt SGEMM} is shown in the inset as a thick
line.

For the smaller basis set, \SpAMM{}($4 \times 4$) is about 1.2--1.4 times
faster compared to \SpAMM{}($16 \times 16$), with this gain diminishing with
increasing tolerance. However, no such performance difference is observed for
the larger basis. \SpAMM{}({\tt SGEMM}) is significantly slower for both error
targets with {\tt MKL} surprisingly outperforming {\tt ACML}, the CPU vendor's
library. \SpAMM{}($4 \times 4$) and \SpAMM{}($16 \times 16$) achieve
approximate linear scaling at $n \approx 600$ in STO-2G and $n \approx 2000$ in
6-31G${}^{**}$. For STO-2G, cross-over between \SpAMM{} and {\tt SGEMM} occurs
at $n \approx 300$ for \SpAMM{}($4 \times 4$), $n \approx 500$ for
\SpAMM{}($16 \times 16$), and $n \approx 900$ for \SpAMM{}({\tt SGEMM}, $\tau
= 4.0 \times 10^{-7}$). For 6-31G${}^{**}$, cross-over occurs at $n \approx
2250$ for \SpAMM{}($4 \times 4$, $\tau = 2.0 \times 10^{-8}$) and \SpAMM{}($16
\times 16$, $\tau = 10^{-7}$), $n \approx 1250$ for \SpAMM{}($4 \times 4$,
$\tau = 2.0 \times 10^{-7}$) and \SpAMM{}($16 \times 16$, $\tau = 10^{-6}$),
and $n \approx 3500$ for \SpAMM{}({\tt SGEMM}, $\tau = 1.0 \times 10^{-7}$).

\begin{figure}
\includegraphics[width=1.0\columnwidth]{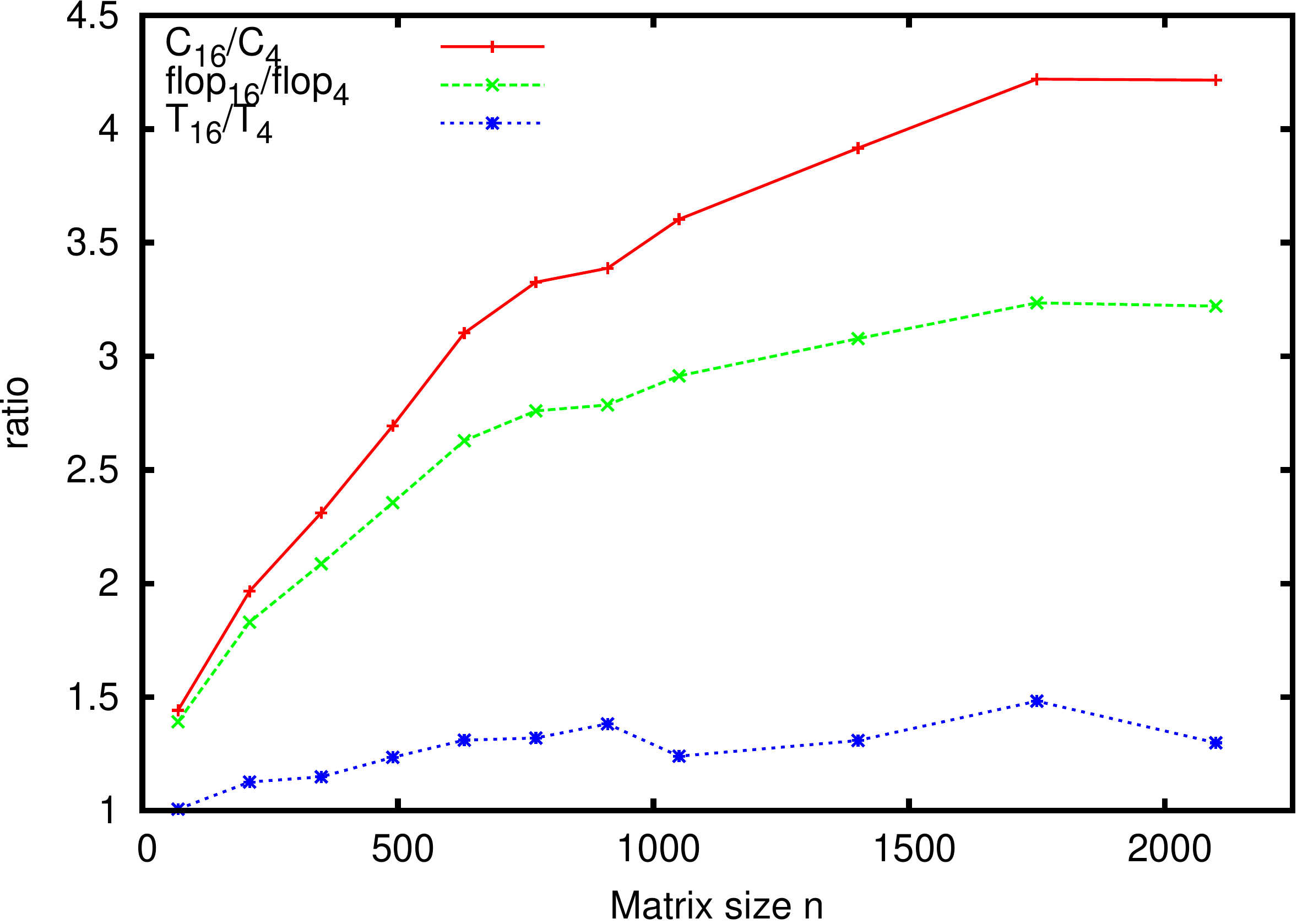}
\caption{\label{fig:water_STO2G_AMD_complexity} Performance ratios of
\SpAMM{}($4 \times 4$, $\tau = 1.0 \times 10^{-7}$) and \SpAMM{}($16 \times
16$, $\tau = 4.0 \times 10^{-7}$) for water clusters of different sizes in
RHF/STO-2G on AMD. Shown are the complexity ratio, $C_{16}/C_{4}$ (the number
of $4 \times 4 \times 4$ matrix products), the flop ratio, $f_{16}/f_{4}$
(measured), and the cycle count ratio, $T_{16}/T_{4}$.}
\end{figure}

\begin{figure}
\includegraphics[width=1.0\columnwidth]{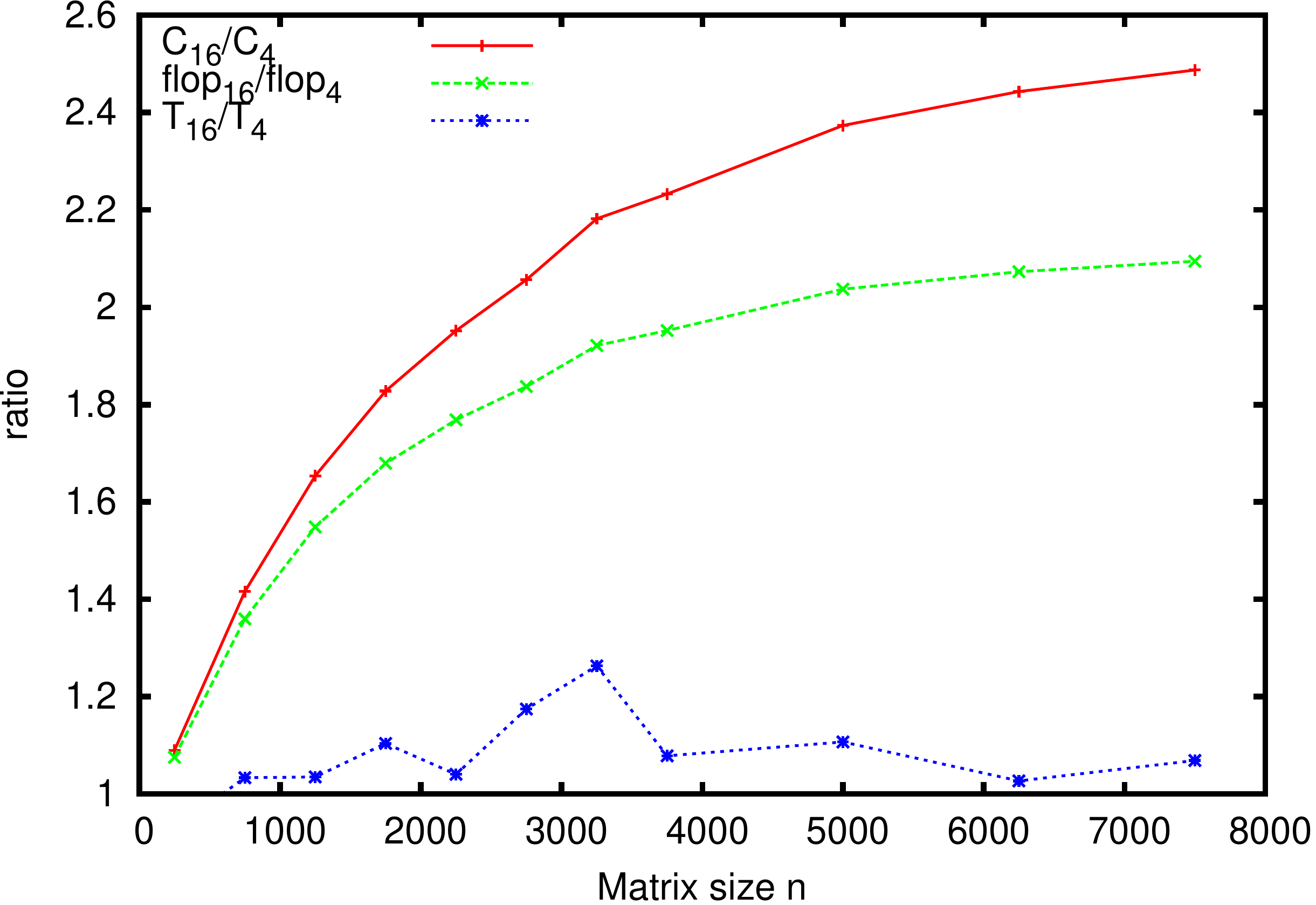}
\caption{\label{fig:water_631gss_AMD_complexity} Performance ratios of
\SpAMM{}($4 \times 4$, $\tau = 2.0 \times 10^{-8}$) and \SpAMM{}($16 \times
16$, $\tau = 10^{-7}$) for water clusters of different sizes in RHF/6-31G** on
AMD.}
\end{figure}

Figures \ref{fig:water_STO2G_AMD_complexity} and
\ref{fig:water_631gss_AMD_complexity} show performance ratios of tests 3 and 4
of Table \ref{table:benchmarks} for a target error of $10^{-6}$. Shown are the
ratios of complexity ($C_{16}/C_{4}$), measured flops
(flop${}_{16}$/flop${}_{4}$), and cycle count ($T_{16}/T_{4}$), where we
define the complexity as the number of $4 \times 4 \times 4$ matrix
multiplications executed, and one $16 \times 16 \times 16$ multiplication
translates into a complexity of 64. For the largest water cluster, the
complexity ratio is about $4.5$ for STO-2G and 2.5 for 6-31G${}^{**}$, the
flop ratio is about 3.7 for STO-2G and 2.1 for 6-31G${}^{**}$, and the cycle
ratio is about 1.3 for STO-2G and 1.1 for 6-31G${}^{**}$. Note that the
difference between complexity and flop ratios is due to the norm product
computations in the conditionals.

\section{Discussion}
\label{sec:discussion}

The error achieved by \SpAMM{}($4 \times 4$, $\tau < 2.0 \times 10^{-8}$)
(Figs.~\ref{fig:water_631gss_RHF_STO_2G_good_error_AMD} and
\ref{fig:water_631gss_RHF_6_31Gss_tight_error_AMD}) is lower than {\tt SGEMM}.
We attribute this behavior to superior order of operation associated with
hierarchical summation of equi-magnitude blocks. A related problem for future
work involves establishing precise relationships between the \SpAMM{}
approximation, matrix decay and application specific measures of error.  In
addition, our implementation can be readily extended to double-precision;
however, our goal in this work is to achieve full single-precision (or better)
and an early onset of linear scaling.

The storage complexity for dense matrices is $\mathcal{O} \left( n^{2}
\right)$ but $\mathcal{O} (n)$ in the sparsified case for matrices exhibiting
approximately exponential decay. The computational complexity for convolution
is $\mathcal{O} \left( n^{2} \log n^{2} \right)$ in the dense and $\mathcal{O}
\left( n \log n \right)$ in the sparsified case due to either octree recursion
or index sorting in the symbolic multiply, see Fig.~\ref{fig:symbolic_cost}.
The computational complexity of the numeric multiply is $\mathcal{O} (n)$ in
either case, when decay is fast enough, or the matrix threshold is large
enough. For the matrices studied here, we find that the numeric multiply
dominates a small but admittedly quadratic component from the symbolic kernel.
Nevertheless we observe an approximately quadratic increase in effective
performance, as shown in Fig.~\ref{fig:water_631gss_sgemm_comparison_AMD},
with an approximately linear increase in cycle count, as shown in
Figs.~\ref{fig:water_STO2G_AMD} and \ref{fig:water_631gss_AMD}.

Comparing the cycle counts shown in Figs.~\ref{fig:water_STO2G_AMD} and
\ref{fig:water_631gss_AMD} between \SpAMM{}($4 \times 4$) and \SpAMM{}($16
\times 16$), \SpAMM{}($4 \times 4$) holds a slight advantage over \SpAMM{}($16
\times 16$) for both basis sets. We attribute differences to variations
between matrix structures; in the STO-2G basis the submatrices of H and O are
small and \SpAMM{}($4 \times 4$) is able to more fully exploit the matrix
structure leading to a bigger performance advantage, while in the
6-31G${}^{**}$ basis, the submatrices of H and O are larger, and \SpAMM{}($16
\times 16$) is more appropriate. We find that the corresponding flop count
ratios, flop${}_{16}$/flop${}_{4}$, (shown in
Figs.~\ref{fig:water_STO2G_AMD_complexity} and
\ref{fig:water_631gss_AMD_complexity} to be significantly larger than the
cycle count ratios, $T_{16}/T_{4}$, clearly indicating the difficulty of
translating the reduced complexity of fine grained algorithms into performance
gains on commodity platforms. It should be noted that both numeric kernels
make full use of SSE, and the issue is not one of inefficient use of floating
point resources.

The significant complexity reduction shown in
Figs.~\ref{fig:water_STO2G_AMD_complexity} and
\ref{fig:water_631gss_AMD_complexity}, as well as the early onset of linear
scaling shown in Figs.~\ref{fig:water_STO2G_AMD} and
\ref{fig:water_631gss_AMD}, underscore the importance of fine granularities.
However, the corresponding irregular access prevents the hardware prefetch
mechanism in the CPU from activating, resulting in a performance loss due to
higher memory latencies, partially offsetting the performance gains due to
complexity reduction (2--3x in this case). So far we have been unsuccessful in
closing this performance gap by using software prefetch statements and $16
\times 16$ blocks, as additional conditionals and the required large distance
between prefetch and data lead to no improvement.  Indeed, it appears that
optimal placement of software prefetch is not only difficult, but that the
performance gains are generally limited \cite{chilimbi2000making,
yang2000push, badawy2004efficacy, huang2009performance}.

Recently much effort has been dedicated to the co-design of improved prefetch
mechanisms, see for instance the ``Scheduled Region Prefetch''
\cite{lin2001reducing, wang2003guided}, ``Flux Cache''
\cite{gaydadjiev2006sad}, ``Traversal Cache'' \cite{stitt2008traversal}, or
the ``Block Prefetch Operation'' \cite{karlsson2000prefetching}.  These works
support and underscore our conclusion that it is the mechanism of memory
access, and not the memory hierarchy \emph{per se} that is the bottle neck for
fine grained irregular algorithms including \SpAMM{}($4 \times 4$), related
$N$-body algorithms, as well as the generic pointer chasing problem.

\section{Conclusions}
\label{sec:conclusions}

We have developed an optimized single-precision implementation of the \SpAMM{}
algorithm using linkless trees and SSE intrinsics in combination with a
hand-coded assembly kernel, and demonstrated early cross-overs with vendor
tuned {\tt SGEMM}s and an early onset of linear scaling for dense quantum
chemical matrices.  For small values of the product space threshold $\tau$, we
find \SpAMM{} exhibits superior error control relative to {\tt SGEMM}.  While
this is unsurprising given the $\mathcal{O} \left( n^{2} \right)$ error
scaling for classical matrix-matrix multiplication relative to an $\mathcal{O}
\left(n \log n \right)$ scaling for recursive multiplication
\cite{Bini:1980:FastMM}, it is worth noting that our fast implementation is
doing something else entirely, namely norm sorting prior to convolution (as
described in Section \ref{sec:symbolic_multiply}). The error analysis of this
approach, as well as for approximate recursive multiplication of SFC clustered
matrices, is an interesting and open problem.  Importantly, the issue has been
refocused on achieving performance and error control that is superior to the
classical multiply, rather than on the errors associated with vector space
truncation and their accumulation in classical multiplication.

In this implementation\footnote{To be released under the GPLv3 at
\url{http://www.freeon.org/}.}, the matrices are stored in contiguous chunks
of size $16 \times 16$, which is unusually small compared to other high
performance matrix multiplies. In addition, the norm condition,
Eq.~(\ref{eq:norm_condition}), is applied to $4 \times 4$ submatrices giving
the implementation fine grained error control with a remaining potential for
2--3x gains in performance on future hardware with improved prefetch.  We
speculate that these gains might be unlocked with hardware modifications such
as the Block Prefetch Operation \cite{karlsson2000prefetching}.  Nevertheless,
relative to na\"{\i}ve implementations of \SpAMM{} using vendor tuned versions
of {\tt SGEMM}, such as those found in Intel's Math Kernel Library ({\tt MKL})
or AMD's Core Math Library ({\tt ACML}), our optimized version is found to be
significantly faster.

Matrix-matrix multiplication is now aligned with the $N$-Body programming
models employed by other quantum chemical solvers
\cite{Challacombe:1996:QCTCb, Challacombe:2000:HiCu, Challacombe:1996:QCTCa},
and more broadly with the generalized $N$-Body framework
\cite{hernquist1988hierarchical, liu1994experiences, Gray:2001:NBody,
Gray:2003:NBody}.  As such, its worth noting that SFC ordering in the vector
and product space, shown in (C) and (D) of Fig.~\ref{fig:SFC}, may also
provide a mechanism for domain decomposition and load balance, corresponding
to methods that have been already proven for other parallel irregular $N$-Body
problems \cite{Warren:1992:HOT, Warren:1995:HOTb, Aluru:1997:SFC,
Campbell:2003:SFC, Devine:2005:SFC}.  Likewise, overdecomposition of the
recursive, three-dimensional {\tt SpAMM} task space may provide significant
opportunities for parallelism through advanced middleware such as {\tt
charm++} \cite{CharmppPPWCPP96} or {\tt tascel} \cite{lifflander2012work,
ma2012data}, as well as through task parallel programming models supported by
{\tt Cilk Plus} \cite{cilkplus} or version 3.0 of the {\tt OpenMP} standard
\cite{openmp}.  The low level SIMD optimizations based on linkless trees
presented here also may be extended to evolving vector lengths through
compiler auto-vectorization and SIMD language extensions such as those
provided by the {\tt ispc} compiler \cite{ispc}, with applicability to AVX,
Intel Xeon Phi architectures and GPUs.

\section*{Acknowledgements}

This work was supported by the U.~S.~ Department of Energy under Contract
No.~DE-AC52-06NA25396 and LDRD-ER grant 20110230ER. The friendly and
stimulating atmosphere of the Ten Bar Caf\'{e} is gratefully acknowledged. NB
gratefully acknowledges the helpful discussions on {\tt irc.freenode.net
\#asm} and {\tt linux-assembly@vger.kernel.org}. NB and MC gratefully
acknowledge helpful discussions with Michele Benzi, Fred Gustavson, and David
Wise. Released under LA-UR 11-06091.

\appendix

\section{Effective Performance}
\label{sec:effective_performance}
Typically the performance of a solver implementation is measured by the number
of CPU cycles (seconds of CPU time), or the number of floating point
operations (flop) per cycle (second). While the former directly measures
``time to solution'', the latter gives insight into the efficiency of the
implementation executing computation instructions under the constraints of
memory access, register pressure, and other factors. Because of the reduction
in flops through Eq.~(\ref{eq:norm_condition}), the perceived performance of
\SpAMM{} (inverse time to solution) is different than its measured floprate.
For comparisons with {\tt SGEMM}, we therefore model the number of flops by $F
=  m \left[ k(1+2n) - n \right]$ and express the effective performance as the
ratio of $F$ and the number of CPU cycles (seconds).

\bibliographystyle{siam}
\bibliography{spammpack_high_performance_implementation}

\end{document}